\documentclass[12pt]{amsart}
\usepackage{amsmath,amscd,amssymb,amsfonts}
\newcommand{\ssb}{\raise.15ex\h{${\scriptscriptstyle\bullet}$}}
\newcommand{\ssc}{\,\raise.15ex\h{${\scriptstyle\circ}$}\,}
\newcommand{\h}{\hbox}
\newcommand{\q}{\quad}

\newcommand{\bs}{\par\bigskip}
\newcommand{\ms}{\par\medskip}
\newcommand{\sk}{\par\smallskip}
\newcommand{\bsn}{\par\bigskip\par\noindent}
\newcommand{\msn}{\par\medskip\par\noindent}
\newcommand{\skn}{\par\smallskip\par\noindent}
\newcommand{\mopl}{\h{$\bigoplus$}}
\newcommand{\msum}{\h{$\sum$}}
\newcommand{\mcap}{\h{$\bigcap$}}
\newcommand{\A}{{\mathcal A}}
\newcommand{\D}{{\mathcal D}}
\newcommand{\E}{{\mathcal E}}
\newcommand{\K}{{\mathcal K}}
\newcommand{\Lc}{{\mathcal L}}
\newcommand{\M}{{\mathcal M}}
\newcommand{\Hc}{{\mathcal H}}
\newcommand{\I}{{\mathcal I}}
\newcommand{\Oc}{{\mathcal O}}
\newcommand{\C}{{\mathbf C}}
\newcommand{\DD}{{\mathbf D}}
\newcommand{\N}{{\mathbf N}}
\newcommand{\PP}{{\mathbf P}}
\newcommand{\Q}{{\mathbf Q}}
\newcommand{\R}{{\mathbf R}}
\newcommand{\Z}{{\mathbf Z}}
\newcommand{\HS}{{\mathbf H\mathbf S}}
\newcommand{\bo}{{\mathbf 1}}

\newcommand{\Mo}{{}\,\overline{\!M}{}}
\newcommand{\So}{{}\,\overline{\!S}{}}
\newcommand{\Xo}{{}\,\overline{\!X}{}}
\newcommand{\Kt}{\widetilde{K}}
\newcommand{\Mt}{\widetilde{M}}
\newcommand{\MMt}{\widetilde{\mathcal M}}
\newcommand{\St}{\widetilde{S}}
\newcommand{\Xt}{\widetilde{X}}
\newcommand{\Zt}{\widetilde{Z}}

\newcommand{\mh}{{}^{\mathfrak m}}
\newcommand{\dd}{\partial}
\newcommand{\ddd}{{\rm d}}
\newcommand{\al}{\alpha}
\newcommand{\la}{\lambda}
\newcommand{\Dec}{{\rm Dec}}
\newcommand{\Om}{\Omega}
\newcommand{\DR}{{\rm DR}}
\newcommand{\Gr}{{\rm Gr}}
\newcommand{\MF}{{\rm MF}}
\newcommand{\MH}{{\rm MH}}
\newcommand{\MHM}{{\rm MHM}}
\newcommand{\MHW}{{\rm MHW}}
\newcommand{\into}{\hookrightarrow}
\newcommand{\simto}{\buildrel{\sim}\over\longrightarrow}
\newcommand{\ges}{\geqslant}
\newcommand{\les}{\leqslant}
\newcommand{\bl}{\bigl}
\newcommand{\br}{\bigr}
\newcommand{\1}{\hskip1pt}
\begin{document}
\title{A young person's guide to mixed Hodge modules}
\author[M. Saito]{Morihiko Saito}
\address{RIMS Kyoto University, Kyoto 606-8502 Japan}
\begin{abstract}
We give a rather informal introduction to the theory of mixed Hodge modules for young mathematicians.
\end{abstract}
\maketitle
\centerline{\bf Introduction}
\bsn
The theory of mixed Hodge modules (\cite{mhp}, \cite{mhm}) was originally constructed as a Hodge-theoretic analogue of the theory of $\ell$-adic mixed perverse sheaves (\cite{th1}, \cite{weil}, \cite{BBD}) and also as an extension of Deligne's mixed Hodge theory (\cite{th2}, \cite{th3}) including the theory of degenerations of variations of pure or mixed Hodge structures (see \cite{Sch}, \cite{St}, \cite{CK}, \cite{CKS1}, \cite{StZ}, \cite{Kadm}, etc.)
The main point is the ``stability" by the direct images and the pull-backs under morphisms of complex algebraic varieties and also by the dual and the nearby and vanishing cycle functors.
Here ``stability" means more precisely the ``stability theorems" saying that there are {\it canonically defined functors} between the derived categories.
These stability theorems become quite useful by combining them with the fundamental theorem of mixed Hodge modules asserting that any admissible variations of mixed Hodge structure on smooth complex algebraic varieties in the sense of \cite{StZ}, \cite{Kadm} are mixed Hodge modules \cite{mhm}; in particular, mixed Hodge modules on a point are naturally identified with graded-polarizable mixed $\Q$-Hodge structures in the sense of Deligne \cite{th2}.
\sk
Technically the {\it strictness} of the Hodge filtration $F$ on the underlying complexes of $\D$-modules is very important in the theory of mixed Hodge modules. Here $\D$-modules are indispensable for the generalization of Deligne's theory in the absolute case (\cite{th2}, \cite{th3}) to the {\it relative} case.
For instance, this includes the assertion that any morphism of mixed Hodge modules is {\it bistrict} for the Hodge and weight filtrations $F,W$, generalizing the case of mixed Hodge structures by Deligne \cite{th2}.
$\D$-modules are also essential for the construction of a relative version of $\Dec\,W$ in the proof of the stability theorem of mixed Hodge modules by the direct images under projective morphisms extending Deligne's argument in the absolute case, where we have {\it bistrict} complexes of $\D$-modules with filtrations $(F,\Dec\,W)$ under the direct images by projective morphisms, see \cite[Proposition 2.15]{mhm}.
\sk
In order to understand the theory of mixed Hodge modules, the general theory of $\D$-modules does not seem to be absolutely indispensable.
In fact, $\D$-modules appearing in the theory of mixed Hodge modules are rather special ones, and we have to deal with them always as {\it $F$-filtered $\D$-modules}, where a slightly different kind of argument is usually required.
For instance, although the {\it regularity} of $\D$-modules is quite useful for the construction of the direct images by affine open immersions, this can be reduced essentially to the normal crossing case by using Beilinson's functor together with the stability theorem under the direct images by projective morphisms (see \cite{def}), and in the latter case, a generalization of the Deligne canonical extensions \cite{eq} to the case of $\D$-modules with normal crossing singular supports is sufficient.
Also the pull-backs of mixed Hodge modules under closed immersions are constructed by using nearby and vanishing cycle functors, which is entirely different from the usual construction of pull-backs of $\D$-modules, see \cite[Section 4.4]{mhm}.
\sk 
Note finally that we need only Zucker's Hodge theory in the curve case \cite{Zu} for the proof of the stability theorem of mixed Hodge  modules under the direct images by proper morphisms, and classical Hodge theory is not used, see (2.5) below.
\sk
I thank the referee for useful comments.
The author is partially supported by Kakenhi 15K04816.
\sk
In Section~1 we explain the main properties of pure Hodge modules.
In Section~2 we give an inductive definition of pure Hodge modules, and explain an outline of proofs of Theorems~(1.3) and (1.4).
In Section~3 we explain a simplified definition of mixed Hodge modules following \cite{def}.
In Appendices we explain some basics of hypersheaves, $\D$-modules, and compatible filtrations.
\msn
{\bf Conventions 1.} We assume that algebraic varieties in this paper are always defined over $\C$ and are reduced (but not necessarily irreducible).
More precisely, a variety means a separated reduced scheme of finite type over $\C$, but we consider only its {\it closed points}, that is, its {\it $\C$-valued points}. So it is close to a variety in the sense of Serre (except that reducible varieties are allowed here). We also assume that a variety is always quasi-projective, or more generally, globally embeddable into a smooth variety (where morphisms of varieties are assumed quasi-projective) in order to simplify some arguments. The reader may assume that all the varieties in this paper are reduced quasi-projective complex algebraic varieties.
\skn
{\bf 2.} We use {\it analytic} sheaves on complex algebraic varieties; in particular, any $\D$-modules are analytic $\D$-modules.
(These are suitable for calculations using local coordinates.)
For the underlying filtered $\D$-module $(M,F)$ of a mixed Hodge module on $X$, one can pass to the corresponding {\it algebraic} filtered $\D$-module by applying GAGA to each $F_pM$ after taking an extension of the mixed Hodge module over a compactification of $X$.
\skn
{\bf 3.} In this paper, perverse sheaves (which do not seem appropriate at least for book titles, see \cite{Di}, \cite{KS}) are mainly called hypersheaves by an analogy with hypercohomology versus cohomology (although the word ``sheaf" may not be of Greek origin). The abelian category of hypersheaves on $X$ are denoted by $\HS(X,A)$ where $A$ is a subfield of $\C$.
Hypersheaves are not sheaves in the usual sense, but they behave like sheaves in some sense; for instance, they can be defined locally provided that gluing data satisfying some compatibility condition are also given.
\bs\bs
\vbox{\centerline{\bf 1. Main properties of pure Hodge modules}
\bsn
In this section we explain the main properties of pure Hodge modules.}
\msn
{\bf 1.1.~Filtered $\D$-modules with $\Q$-structure.}
A pure Hodge module $\M$ on a smooth complex algebraic variety $X$ of dimension $d_X$ is basically a coherent left $\D_X$-module $M$ endowed with the Hodge filtration $F$ such that $(M,F)$ is a filtered left $\D_X$-module with $F_pM$ coherent over $\Oc_X$.
The filtration $F$ on $\D_X$ is by the order of differential operators, and we have the filtered $\D$-module condition
$$(F_p\D_X)\,F_qM\subset F_{p+q}M\q\q(p\in\N,\,q\in\Z),
\leqno(1.1.1)$$
(which is equivalent to {\it Griffiths transversality} in the case of variations of Hodge structure), and the equality holds in (1.1.1) for $q\gg0$.
Moreover $M$ has a $\Q$-{\it structure} given by an isomorphism
$$\alpha_{\M}:\DR_X(M)\cong\C\otimes_{\Q}K\q\h{in}\,\,\,\HS(X,\C),
\leqno(1.1.2)$$
with $K\in\HS(X,\Q)$ (see (A.1) below). Here $\DR_X(M)$ is the de Rham complex of the $\D_X$-module $M$ viewed as a quasi-coherent $\Oc_X$-module with an integrable connection:
$$\DR_X(M):=\bl[M\to\Om_X^1\otimes_{\Oc_X}M\to\cdots\to\Om_X^{d_X}\otimes_{\Oc_X}M\br],
\leqno(1.1.3)$$
with last term put at degree 0. In (1.1.2) we {\it assume}
$$\DR_X(M)\in\HS(X,\C),
\leqno(1.1.4)$$
and moreover (1.1.4) holds with $M$ replaced by any subquotients of $M$ as coherent $\D$-modules. (More precisely, these properties follow from the condition that $M$ is {\it regular holonomic}, see (B.5.2) below. We effectively assume the last condition in the definition of Hodge modules, see Remark~(ii) after Theorem~(1.3).)
In this paper we also assume
$$\h{$M,K$ are quasi-unipotent (see (B.6) below).}
\leqno(1.1.5)$$
\sk
We will denote by $\MF(X,\Q)$ the category of
$$\M=\bl((M,F),K,\alpha_{\M}\br)$$
satisfying the above conditions.
(Sometimes $\alpha_{\M}$ will be omitted to simplify the notation.)
We also use the notation
$${\rm rat}(\M):=K.
\leqno(1.1.6)$$
The category $\MH(X,w)$ of pure Hodge modules of weight $w$ on $X$ will be defined as a full subcategory of $\MF(X,\Q)$.
\msn
{\bf 1.2.~Strict support decomposition.}
The first condition for pure Hodge modules is the {\it decomposition by strict support}
$$\M=\mopl_{Z\subset X}\,\M_Z\q\h{in}\,\,\,\MF(X,\Q),
\leqno(1.2.1)$$
where $Z$ runs over irreducible closed subvarieties of $X$, and
$$\M_Z=\bl((M_Z,F),K_Z,\al_{\M_Z}\br)$$
has {\it strict support} $Z$ (that is, its support is $Z$, and it has no nontrivial sub nor quotient object supported on a proper subvariety of $Z$).
More precisely, we assume that the last condition is satisfied for both $M_Z$ and $K_Z$.
Note that the condition for $K_Z$ is equivalent to that $K_Z$ is an intersection complex with local system coefficients, see \cite{BBD}.
We then get
$${\rm Hom}(\M_Z,\M_{Z'})=0\q\h{if}\,\,\,Z\ne Z'.
\leqno(1.2.2)$$
In fact, this holds with $\M_Z,\M_{Z'}$ replaced by $M_Z,M_{Z'}$ or by $K_Z,K_{Z'}$.
\sk
In (2.2) below, we will define the full subcategory
$$\MH_Z(X,w)\subset\MF(X,\Q)$$
consisting of {\it pure Hodge modules of weight $w$ with strict support $Z$} by increasing induction on $\dim Z$, and put
$$\MH(X,w):=\mopl_{Z\subset X}\,\MH_Z(X,w)\subset\MF(X,\Q),
\leqno(1.2.3)$$
where the direct sum over closed irreducible subvarieties $Z$ of $X$ is justified by (1.2.1--2).
\sk
This full subcategory $\MH_Z(X,w)\subset\MF(X,\Q)$ can be defined effectively by the following {\it fundamental theorem of pure Hodge modules}, which may be viewed as a {\it working definition} of $\MH_Z(X,w)$.
\msn
{\bf Theorem~1.3} (\cite[Theorem~3.21]{mhm}). {\it For any closed irreducible subvariety $Z\subset X$, the restriction to sufficiently small open subvarieties of $Z$ induces an equivalence of categories
$$\MH_Z(X,w)\simto {\rm VHS}_{\rm gen}(Z,w-\dim Z)^p,
\leqno(1.3.1)$$
where the right-hand side is the category of polarizable variations of pure Hodge structure of weight $w-\dim Z$ defined on smooth dense open subvarieties $U$ of $Z$. $($More precisely, we take the inductive limit over $U\subset Z.)$
Moreover, $(1.3.1)$ induces a one-to-one correspondence between polarizations of $\M\in\MH_Z(X,w)$ $($see $(2.2.2)$ below$)$ and those of the corresponding generic variation of Hodge structure.}
\msn
{\bf Remarks.} (i) The equivalence of categories (1.3.1) means that any pure Hodge module with strict support $Z$ is generically a polarizable variation of pure Hodge structure, and conversely any polarizable variation of pure Hodge structure defined on a smooth dense Zariski-open subset $U\subset Z$ can be extended uniquely to a pure Hodge module with strict support $Z$.
\sk
(ii) By using the stability theorem of pure Hodge modules under the direct images by projective morphisms (see Theorem~(1.4) below), the proof of Theorem~(1.3) can be reduced to the case where $Z=X$ and a variation of Hodge structure is defined on the complement $U$ of a divisor with normal crossings. In this case we use the original definition of pure Hodge modules in \cite{mhp} given by induction on the dimension of the support of $M$ and using the nearby and vanishing cycle functors.
\sk
In the case $D:=X\setminus U$ is a divisor with normal crossings, the above extension of a variation of Hodge structure on $U$ to a Hodge module on $X$ is rather easy to describe by using Deligne's canonical extension \cite{eq}. In fact, we have the following explicit formula (see \cite[(3.10.12)]{mhm}:
$$F_pM=\msum_{i\ges 0}\,F_i\1\D_X\1\bl(j_*F_{p-i}\1\Lc\cap\widehat{\Lc}^{\,>-1}\br)\subset\widehat{\Lc}^{\,>-1}(*D),
\leqno(1.3.2)$$
with $j:U\into X$ is the inclusion. Here $\Lc$ is the locally free $\Oc_U$-module underlying the variation of Hodge structure with $F$ the Hodge filtration, $\widehat{\Lc}^{\,>-1}$ is the Deligne extension of $\Lc$ over $X$ such that the eigenvalues of the residues of the connection are contained in $(-1,0]$, and the last term of (1.3.2) is the Deligne meromorphic extension, see \cite{eq}. (Note that a decreasing filtration $F$ is identified with an increasing filtration by setting $F_p:=F^{-p}$.)
Taking the union over $p\in\Z$, we have
$$M=\D_X\,\widehat{\Lc}^{\,>-1}\subset\widehat{\Lc}^{\,>-1}(*D).
\leqno(1.3.3)$$
\sk
The underlying $\Q$-local system $L$ of a polarizable variation of Hodge structure on $U$ is canonically extended over $X$ as an intersection complex (see \cite{BBD}) where $L$ must be shifted by $\dim Z$.
This can be done without assuming that $D=X\setminus U$ is a divisor with normal crossings on a smooth variety.
\ms
The second main theorem in the theory of pure Hodge modules is the stability theorem of pure Hodge modules under the direct images by projective morphisms.
\msn
{\bf Theorem~1.4} (\cite[Theorem~5.3.1]{mhp}). {\it Let $f:X\to Y$ be a projective morphism of smooth complex algebraic varieties, and $\M=((M,F),K,\alpha_{\M})\in\MH_Z(X,w)$. Let $\ell$ be the first Chern class of an $f$-ample line bundle. Then the direct image $f_*^{\D}(M,F)$ as a filtered $\D$-module $($see {\rm (B.3)} below$)$ is strict, and we have
$$\Hc^if_*\M:=(\Hc^if_*^{\D}(M,F),\mh\Hc^if_*K,\mh\Hc^if_*\alpha_{\M})\in\MH(Y,w+i)\q(i\in\Z),
\leqno(1.4.1)$$
together with the isomorphisms
$$\ell^i:\Hc^{-i}f_*\M\simto\Hc^if_*\M(i)\q(i>0),
\leqno(1.4.2)$$
where $(i)$ denotes the Tate twist shifting the filtration $F$ by $i$, see a remark after $(2.1.10)$ below.
\sk
Moreover, if $S:K\otimes K\to\DD_X(-w)$ is a polarization of $\M$ {\rm(}see $(2.2.2)$ below$)$, then a polarization of the $\ell$-primitive part
$${}^P\Hc^{-i}f_*\M:={\rm Ker}\,\ell^{\,i+1}\subset\Hc^{-i}f_*\M\q(i\ges 0)$$
is given by the restriction to the $\ell$-primitive part of the induced pairing}
$$(-1)^{i(i-1)/2}\,\mh\Hc f_*S\ssc(id\otimes\ell^{\,i}):\mh\Hc^{-i}f_*K\otimes\mh\Hc^{-i}f_*K\to\DD_Y(i-w).
\leqno(1.4.3)$$
\msn
{\bf Remarks.} (i) The action of $\ell$ can be defined by using $C^{\infty}$ forms and the Dolbeault resolution to get a filtered complex which is filtered quasi-isomorphic to the relative de Rham complex
$$\DR_{X\times Y/Y}\bl((i_f)_*^{\D}(M,F)\br),$$
which is used in the definition of the direct image in (B.3) below. Note that the first Chern class is represented by a closed 2-form of type $(1,1)$ on $X$, see also \cite[Lemma 5.3.2]{mhp}.
(It is also possible to define the action of $\ell$ by using the restriction to a sufficiently general hyperplane section of $(M,F)$.)
\sk
(ii) We use Deligne's sign convention of a polarization $S$ of a Hodge structure $(H_{\C},F),H_{\Q})$, see \cite[Definition 2.1.15]{th2}; that is,
$$S(v,C\,\overline{\!v})>0\q(v\in H_{\C}\setminus\{0\}),
\leqno(1.4.4)$$
where $C$ is the Weil operator defined by $i^{p-q}$ on $H^{p,q}_{\C}$, and the Tate twist $(2\pi i)^w$ is omitted to simplify the notation.
(This sign convention seems to be theoretically natural if one considers the action of the Weil restriction of ${\mathbf G}_{\mathbf m}$, see \cite[Section 2.1]{th1}.)
\sk
Recall that the usual sign convention is
$$S(C\1v,\,\overline{\!v})>0\q(v\in H_{\C}\setminus\{0\}),
\leqno(1.4.5)$$
and the difference is given by the multiplication by
$$(-1)^w,
\leqno(1.4.6)$$
where $w$ is the weight of the Hodge structure.
\sk
If we use the usual sign convention (1.4.5) instead of (1.4.4), then the difference (1.4.6) implies a considerable change of the formula (1.4.3) of Theorem~(1.4). In fact, we have to change the sign for each direct factor with strict support $Z$ {\it depending on} $d_Z:=\dim Z$, since the difference of sign (1.4.6) depends on the {\it pointwise weight}
$$w-i-d_Z
\leqno(1.4.7)$$
of the generic variations of Hodge structure of each direct factor with strict support $Z$:
$$(\Hc^{-i}f_*\M)_Z\subset\Hc^{-i}f_*\M.$$
\sk
(iii) For a polarization $S$ of a generic variation of Hodge structure of a pure Hodge module with strict support $Z$, the associated polarization of the Hodge modules is defined by
$$(-1)^{d_Z(d_Z-1)/2}\,S:K|_U\otimes K|_U\to\DD_U(-w),
\leqno(1.4.8)$$
on an open subvariety $U\subset Z$ where the variation of Hodge structure is defined, see \cite[Proposition 5.2.16]{mhp} for the constant coefficient case. This can be extended to a pairing of $K$ by using the theory of intersection complexes \cite{BBD}.
\sk
For a smooth projective variety $X$, it is well-known that
$$\int_X(-1)^{j(j-1)/2}\,i^{p-q}\,v\,\wedge\,\overline{\!v}{}\wedge\omega^{d_X-j}>0,
\leqno(1.4.9)$$
for a nonzero element $v$ in the primitive cohomology ${}^P\!H^j(X,\C)$ of type $(p,q)$, where $\omega$ is a K\"ahler form. One may think that (1.4.9) contradicts (1.4.8) combined with (1.4.3) in Theorem~(1.4) for $f:X\to pt$. In fact, by setting $\xi(k)=k(k-1)/2$ for $k\in\N$, the difference between the sign $(-1)^{\xi(j)}$ in (1.4.9) and the product of $(-1)^{\xi(d_X)}$ in (1.4.8) with $Z=X$ and $(-1)^{\xi(d_X-j)}$ in (1.4.3) with $i=d_X-j$ does not coincide with the difference between the two sign conventions (1.4.4--5) which is equal to $(-1)^j$ by (1.4.6) for $w=j$, since
$$\xi(d_X)-\xi(d_X-j)-\xi(j)-j=jd_X\mod 2.$$
However, the remaining sign $(-1)^{jd_X}$ just comes from the isomorphism
$$\R\Gamma(X,\C_X)[d_X]\otimes_{\C}\R\Gamma(X,\C_X)[d_X]\cong\bl(\R\Gamma(X,\C_X)\otimes_{\C}\R\Gamma(X,\C_X)\br)[2d_X].
\leqno(1.4.10)$$
In fact, if we set $\K^{\ssb}:=\R\Gamma(X,\C_X)$, then $\K^{\ssb}[d_X]$ is identified with $\C_X[d_X]\otimes_{\C}\K^{\ssb}$, and the above isomorphism is given by
$$\C[d_X]\otimes_{\C}\K^{\ssb}\otimes_{\C}\C[d_X]\otimes_{\C}\K^{\ssb}\cong\C[d_X]\otimes_{\C}\C[d_X]\otimes_{\C}\K^{\ssb}\otimes_{\C}\K^{\ssb},$$
where the middle two components are exchanged, but the remaining ones are unchanged (see \cite{coh} for the sign about single complexes associated with $n$-ple complexes).
A similar argument is used in the proof of the anti-commutativity of the last diagram in \cite[Section 5.3.10]{mhp} where $d_X=1$.
\bs\bs
\vbox{\centerline{\bf 2. Outline of proofs of Theorems~(1.3) and (1.4)}
\bsn
In this section we give an inductive definition of pure Hodge modules (see \cite{via}, \cite{mhp}), and explain an outline of proofs of Theorems~(1.3) and (1.4).}
\msn
{\bf 2.1.~Admissibility condition along $g=0$.} Let $X$ be a smooth algebraic variety, and $Z$ be an irreducible closed subvariety of $X$. Let
$$\M=((M,F),K,\alpha_{\M})\in\MF(X,\Q),$$
with strict support $Z$, see (1.2). Let $g$ be a function on $X$, that is, $g\in\Gamma(X,\Oc_X)$. Let $i_g:X\into X\times\C$ be the graph embedding by $g$. Set
$$(\Mt,F):=(i_g)_*^{\D}(M,F),$$
see (B.3) below for $(i_g)_*^{\D}$. (Note that the filtration $F$ is {\it shifted by} $1$, which is the codimension of the embedding.) We have the filtration $V$ on $\Mt$ (see (B.6) below).
\msn
{\bf Definition.} We say that $(M,F)$ is {\it admissible along} $g=0$ (or $g$-{\it admissible} for short) in this paper if the following two conditions are satisfied:
$$t(F_pV^{\al}\Mt)=F_pV^{\al+1}\Mt\q\q(\forall\,\al>0),
\leqno(2.1.1)$$
\vskip-6mm
$$\dd_t(F_p\Gr_V^{\al}\Mt)=F_{p+1}\Gr_V^{\al-1}\Mt\q\q(\forall\,\al<1,\,p\in\Z),
\leqno(2.1.2)$$
see \cite[3.2.1]{mhp}. (These properties were first found in the one-dimensional case, see \cite{hfgs2}.)
\sk
In the case $Z\subset g^{-1}(0)$, $(M,F)$ is $g$-admissible if and only if the following condition is satisfied:
$$g\,F_pM\subset F_{p-1}M\q\q(\forall\,p\in\Z),
\leqno(2.1.3)$$
see \cite[Lemma 3.2.6]{mhp}.
\sk
In the case $Z\not\subset g^{-1}(0)$, let $j:X\times\C^*\into X\times\C$ be the natural inclusion. We have the isomorphisms
$$F_p\Mt=\msum_{i\ges 0}\,\dd_t^i(j_*F_{p-i}\Mt\cap V^{>0}M)\q\q(p\in\Z),
\leqno(2.1.4)$$
if $(M,F)$ is $g$-admissible and moreover the following condition holds:
$$\dd_t:\Gr_V^1(\Mt,F)\to\Gr_V^0(\Mt,F[-1])\,\,\,\h{is strictly surjective,}
\leqno(2.1.5)$$
see \cite[Remark 3.2.3]{mhp}. This is closely related to (1.3.2). (Forgetting $F$, condition (2.1.5) is equivalent to that $M$ has no nontrivial quotient supported in $g^{-1}(0)$, see (B.6.6) below. The strictness of $F$ follows from the properties of Hodge modules as is seen in (2.3.5) below.)
\sk
Assume $(M,F)$ is $g$-admissible. We have the {\it nearby and vanishing cycle functors} $\psi_g$, $\varphi_g$ defined by
$$\psi_g(M,F):=\mopl_{\la\in\C^*_1}\,\psi_{g,\la}(M,F),\q\varphi_g(M,F):=\mopl_{\la\in\C^*_1}\,\varphi_{g,\la}(M,F),
\leqno(2.1.6)$$
\vskip-5mm
$$\aligned\psi_{g,{\mathbf e}(-\al)}(M,F)&:=\Gr_V^{\al}(\Mt,F)\q(\al\in(0,1]),\\\varphi_{g,1}(M,F)&:=\Gr_V^0(\Mt,F[-1]),\endaligned
\leqno(2.1.7)$$
where $\C^*_1:=\{\la\in\C^*\mid|\la|=1\}$, ${\mathbf e}(-\al):=\exp(-2\pi i\al)$, and $\psi_{g,\la}=\varphi_{g,\la}$ ($\la\ne 1$) as in (A.2.8) below.
We have $\varphi_{g,1}(M,F)=(M,F)$ by \cite[Lemma 3.2.6]{mhp} if ${\rm supp}\,M\subset g^{-1}(0)$ and $g\1F_pM\subset F_{p-1}M$ ($p\in\Z$). Note that $F$ is shifted by 1 when the direct image $(i_g)_*^{\D}$ is taken. This is the reason for which $F$ is shifted for $\varphi$ in (2.1.7), and not for $\psi$ (for left $\D$-modules).
\sk
Combining these with (B.6.7) below, we get
$$\aligned\psi_g\M&:=(\psi_g(M,F),\mh\psi_gK,\mh\psi_g\alpha_{\M}),\\\varphi_g\M&:=(\varphi_g(M,F),\mh\varphi_gK,\mh\varphi_g\alpha_{\M})\q\h{in}\,\,\,\,\MF(X,\Q).\endaligned
\leqno(2.1.8)$$
We have similarly $\psi_{g,1}\M$, $\varphi_{g,1}\M$ together with the morphisms
$$\aligned{\rm can}:\psi_{g,1}\M\to\varphi_{g,1}\M,\q&{\rm Var}:\varphi_{g,1}\M\to\psi_{g,1}\M(-1),\\
N:\psi_g\M\to\psi_g\M(-1),\q&N:\varphi_{g,1}\M\to\varphi_{g,1}\M(-1),\endaligned
\leqno(2.1.9)$$
such that the restrictions of $\,{\rm can}$, $\,{\rm Var}$, $N$ to the $\D$-module part $\Gr_V^1\M$, $\Gr_V^0\M$, $\Gr_V^{\al}\M$ ($\al\in[0,1]$) are respectively given by $-\dd_t$, $t$, $-(\dd_tt-\al)$, and we have
$${\rm Var}\ssc{\rm can}=N\,\,\,\,\h{on}\,\,\,\varphi_{g,1}\M,\q{\rm can}\ssc{\rm Var}=N\,\,\,\,\h{on}\,\,\,\varphi_{g,1}\M.
\leqno(2.1.10)$$
Here the Tate twist $(k)$ for $k\in\Z$ in general is essentially the shift of the filtration $F$ by $[k]$. For the $\Q$-coefficient part, it is defined by the tensor product of $\Q(k):=(2\pi i)^k\Q\subset\C$ over $\Q$, see \cite[Definition 2.1.13]{th2}. (Similarly with $\Q$ replaced by any subfield $A\subset\C$)
\sk
By \cite[Lemma 5.1.4]{mhp} we see that the strict support decomposition (1.2.1) holds if for any $g\in\Gamma(U,\Oc_U)$ with $U$ an open subvariety of $X$, $(M,F)|_U$ is $g$-admissible and moreover
$$\varphi_{g,1}\M|_U={\rm Im}\,{\rm can}\oplus{\rm Ker}\,{\rm Var}\q\h{in}\q\MF(U,\Q),
\leqno(2.1.11)$$
\msn
{\bf 2.2.~Inductive definition of pure Hodge modules.} For a smooth complex algebraic variety $X$ and an irreducible closed subvariety $Z\subset X$,
we define the full subcategory $\MH_Z(X,w)\subset\MF(X,\Q)$ by increasing induction on $d_Z:=\dim Z$ as follows (see \cite{via}, \cite{mhp}):
\msn
{\bf Case 1.} If $Z$ is a point $x\in X$, then we have an equivalent of categories
$$\aligned&(i_x)_*:{\rm HS}(w)^p\simto\MH_{\{x\}}(X,w)\\ \h{with}\q\q&(i_x)_*\bl((H_{\C},F),H_{\Q}\br)=\bl((i_x)_*^{\D}(H_{\C},F),(i_x)_*H_{\Q}\br),\endaligned
\leqno(2.2.1)$$
where $i_x:\{x\}\into X$ denotes the canonical inclusion, and ${\rm HS}(w)^p$ denotes the category of polarizable $\Q$-Hodge structures of weight $w$ (see \cite{th2}). The latter is naturally identified with a full subcategory of $\MF(\{x\},\Q)$ (by setting $F_p=F^{-p}$ as usual).
\msn
{\bf Case 2.} If $d_Z>0$, then $\M=((M,F),K)\in\MF(X,\Q)$ with strict support $Z$ belongs to $\MH_Z(X,w)$ if there is a perfect pairing (see (A.3) below)
$$S:K\otimes_{\Q}K\to\DD_X(-w)=\Q_X(d_X-w)[2d_X],
\leqno(2.2.2)$$
which is called a {\it polarization} of $\M$, and the following two conditions are satisfied:
\msn
(i) The pairing $S$ is compatible with the Hodge filtration $F$ in the following sense:
\skn
There is an isomorphism of filtered $\D$-modules
$$\DD(M,F)=(M,F)(w),$$
which corresponds (by using (B.4.6) below) to an isomorphism defined over $\Q$:
$$\DD(K)=K(w),$$
and the latter is identified with the perfect pairing $S$ via (A.3.1) below.
\msn
(ii) For any Zariski-open subset $U\subset X$ and $g\in\Gamma(U,\Oc_U)$, the restriction of $(M,F)$ to $U$ is $g$-admissible, and moreover, in the case $Z\cap U\not\subset g^{-1}(0)$, we have
$$\Gr_k^W\psi_g\M|_U,\,\Gr_k^W\varphi_{g,1}\M|_U\in\MH_{<d_Z}(U,w),
\leqno(2.2.3)$$
$$\h{$\mh\psi_gS\ssc(id\otimes N^i)$ gives a polarization of ${}^P\Gr^W_{w-1+i}\psi_g\M|_U\,\,\,(i\ges 0)$}.
\leqno(2.2.4)$$
\sk
The last term of (2.2.3) is the direct sum of $\MH_{Z'}(U,w)$ with $Z'$ running over closed irreducible subvarieties of $U$ with $d_{Z'}<d_Z$. The weight filtrations $W$ on $\psi_g\M|_U$, $\varphi_{g,1}\1\M|_U$ are the {\it monodromy filtration} associated with the action of $N:=(2\pi i)^{-1}\log T_u$ (see \cite{weil}) which are shifted by $w-1$ and $w$ respectively. This means that $W$ on $\psi_g\M|_U$ {\it with filtration $F$ forgotten} is uniquely determined by the following conditions:
$$\aligned N(W_i\psi_g\M|_U)&\subset(W_{i-2}\psi_g\M|_U)(-1)\q(i\in\Z),\\ N^i:\Gr_{w-1+i}^W\psi_g\M|_U&\simto(\Gr^W_{w-1-i}\psi_g\M|_U)(-i)\q(i\in\N).\endaligned
\leqno(2.2.5)$$
Here the Tate twists may be neglected since $F$ is forgotten in (2.2.5).
(However, the last isomorphism of (2.2.5) is strictly compatible with $F$ by (2.3.3) below if $\psi_g\M|_U$ belongs to $\MHW(U)$ in (2.3) below.)
A similar assertion holds for $\varphi_{g,1}\M|_U$ with $w-1$ replaced by $w$.
\sk
In (2.2.4), the primitive part ${}^P\Gr^W_{w-1+i}\psi_g\M|_U$ is defined by
$${}^P\Gr^W_{w-1+i}\psi_g\M|_U:={\rm Ker}\,N^{i+1}\subset\Gr^W_{w-1+i}\psi_g\M|_U,
\leqno(2.2.6)$$
by using the induced filtration $F$ on the kernel.
For $\psi_gS$, see (A.3.2) below. Note that the condition for a polarization $S$ is also by induction on $d_Z$. For each direct factor of ${}^P\Gr^W_{w-1+i}\psi_g\M|_U$ with $0$-dimensional strict support, we assume that $\psi_gS\ssc(id\otimes N^i)$ induces a polarization of $\Q$-Hodge structure in the sense of \cite{th2} where the place of the Weil operator is different from the usual one as is noted in Remark (ii) after Theorem (1.4).
\msn
{\bf 2.3. Some properties of pure Hodge modules.}
Let $\MH(X,w)$ be as in (1.2.3). Let
$$\MHW(X)$$
be the category of {\it weakly mixed Hodge modules\,} consisting of $(\M,W)$ with $\M\in\MF(X,\Q)$ and $W$ a finite increasing filtration of $\M$, which satisfy
$$\Gr_w^W\M\in\MH(X,w)\q(\forall\,w\in\Z).
\leqno(2.3.1)$$
We have by definition the stability of pure Hodge modules by the nearby and vanishing cycle functors:
$$\psi_g\M|_U,\,\varphi_{g,1}\M|_U\in\MHW(U).
\leqno(2.3.2)$$
It is easy to show the following (see \cite[Proposition 5.1.14]{mhp}):
\msn
(2.3.3)\,\,\,\, $\MHW(X)$ and $\MH(X,w)$ ($w\in\Z$) are abelian categories such that
\sk\q\q\q
any morphisms are strictly compatible with $(F,W)$ and $F$ respectively.
\ms
This assertion is proved by using
$${\rm Hom}(\M,\M')=0\,\,\,\h{if}\,\,\,\M\in\MH(X,w),\,\M'\in\MH(X,w')\,\,\,\h{with}\,\,\,\,w>w'.
\leqno(2.3.4)$$
This is reduced to \cite{th2} by using the assertion that any $\M\in\MH_Z(X,w)$ is generically a variation of Hodge structure of weight $w-d_Z$.
(The latter is an easy part of Theorem~(1.3).)
\sk
These assertions hold without assuming polarizability (see \cite[Section 5.1]{mhp}), and imply
$${\rm can}:\psi_{g,1}\M\to\varphi_{g,1}\M\,\,\,\,\h{is strictly surjective for $(F,W)$},
\leqno(2.3.5)$$
This assures the strict surjectivity in (2.1.5). It also gives a reason for which condition~(2.2.4) is imposed only for $\psi$.
\ms
We prove Theorems~(1.3) and (1.4) by induction on $\dim Z$ using the following rather technical key theorem:
\msn
{\bf Theorem~2.4.} {\it Let $f:X\to Y$ be as in Theorem~$(1.4)$. Let $g\in\Gamma(Y,\Oc_Y)$. Put $h:=fg$. Let $\M=((M,F),K)\in\MF(X,\Q)$ with strict support $Z\not\subset h^{-1}(0)$. Let $S:K\otimes K\to\DD_X(-w)$ be a perfect pairing compatible with the filtration $F$ as in condition~{\rm (ii)} in Case $2$ of $(2.2)$.
Assume that $(M,F)$ is $h$-admissible, we have
$$\Gr_i^W\psi_h\M,\,\Gr_i^W\varphi_{h,1}\M\in\MH(X,i)\q(\forall\,i\in\Z),
\leqno(2.4.1)$$
with $W$ as in $(2.2.5)$, and the conclusions of Theorem~$(1.4)$ are satisfied for the $N$-primitive part
$${}^{P_N}\Gr_{w-1+i}^W\psi_h\M\,\,\,\,\,\h{with polarization}\,\,\,\,\,\psi_hS\ssc(id\otimes N^i)\q(i\ges 0).
\leqno(2.4.2)$$
Then
\skn
{\rm (i)} The filtered direct image $\,f_*^{\D}(M,F)$ is strict on a sufficiently small neighborhood of $g^{-1}(0)$ $($in the classical topology$)$, and the $\Hc^if_*^{\D}(M,F)$ are $g$-admissible.
\skn
{\rm (ii)} The shifted direct image filtration $f_*^{\D}W[j]$ induces the monodromy filtration shifted by $w+j-1$ on
$$\psi_g\Hc^jf_*^{\D}\M=\Hc^jf_*^{\D}\psi_h\M\q(\forall\,j\in\Z),
\leqno(2.4.3)$$
which is denoted by $W$ so that
$$\Gr_i^W\psi_g\Hc^jf_*^{\D}\M\in\MH(Y,i)\q(\forall\,i\in\Z).
\leqno(2.4.4)$$
\skn
{\rm (iii)} We have isomorphisms on a sufficiently small neighborhood of $g^{-1}(0)$:
$$\ell^j:\Hc^{-j}f_*^{\D}\M\simto(\Hc^jf_*^{\D}\M)(j)\q(\forall\,j\ges 0).
\leqno(2.4.5)$$
\skn
{\rm (iv)} On the bi-primitive part ${}^{P_{\ell}}{}^{P_N}\Gr_{w-1-j+i}^W\psi_g\Hc^{-j}f_*^{\D}\M$ defined by
$${\rm Ker}\,\ell^{j+1}\cap{\rm Ker}\,N^{i+1}\subset\Gr_{w-1-j+i}^W\psi_g\Hc^{-j}f_*^{\D}\M,
\leqno(2.4.6)$$
we have a polarization of Hodge module given by the induced pairing}
$$(-1)^{j(j-1)/2}\1\Gr^W\psi_g\Hc f_*S\ssc(id\otimes N^i\ell^j)\q(\forall\,i,j\ges 0).
\leqno(2.4.7)$$
\ms
We first explain how the above theorem is used in the proofs of Theorems~(1.3) and (1.4).
\msn
{\bf 2.5.~Outline of proofs of Theorems~(1.3) and (1.4).} We show the assertions by increasing induction on the dimension of the strict support $Z$.
The order of the induction is rather complicated as is explained below:
\sk
Assume Theorems~(1.3) and (1.4) are proved for $\dim Z<d$. Then Theorem~(1.4) for $\dim Z=d$ with $f(Z)\ne pt$ follows from Theorem~(2.4) where the decomposition by strict support follows from (2.1.11) (which is satisfied by using \cite[Corollary 4.2.4]{mhp}). Using this, we can reduce the proof of Theorem~(1.3) to the normal crossing case where the singular locus of $M$ is a divisor with normal crossings. Here the filtration $F$ can be defined by (1.3.2), and we have to show that the conditions for Hodge modules are satisfied for any locally defined functions $g$. This can be further reduced by using Theorem~(2.4) (but not Theorem~(1.4)) to the case where the union of $g^{-1}(0)$ and the singular locus of $M$ is a divisor with normal crossings. Here we can calculate explicitly the nearby and vanishing cycle functors together with the induced pairing (although these are rather complicated), see \cite{mhm} for details.
\sk
Now we have to prove Theorem~(1.4) in the case $\dim Z=d$ and $f(Z)=pt$.
Here we may assume that $X=\PP^n$, $Y=pt$. Let $\pi:\Xt\to X$ be the blow-up along the intersection $Z$ of two sufficiently general hyperplanes of $\PP^n$. By Theorem~(1.3) for $\dim Z=d$, there is
$$\MMt=((\Mt,F),\Kt)\in\MH(\Xt,w)\q\h{with}\q\MMt|_{\Xt\setminus\pi^{-1}(Z)}=\M|_{X\setminus Z}.$$
We have a polarization $\St$ of $\MMt$ extending the restriction of a polarization $S$ of $\M$ to the complement of $Z\subset X$.
\sk
By Theorem~(1.4) for $f(Z)\ne 0$, we see that $\M$ is a direct factor of $\Hc^0\pi_*\M$. Moreover its complement is isomorphic to $\M_Z(-1)$ since $Z$ is sufficiently general, where $\M_Z$ is the noncharacteristic restriction of $\M$ to $Z$; more precisely,
$$\M_Z=((M_Z,F),K_Z[-2]),$$
with $(M_Z,F)$ and $K_Z$ respectively the noncharacteristic restrictions of $(M,F)$ and $K$ to $Z$.
\sk
So we have the direct sum decomposition
$$\Hc^0\pi_*\MMt=\M\oplus\M_Z(-1),
\leqno(2.5.1)$$
with $\Hc^j\pi_*\MMt=0$ for $j\ne 0$. Here $\pi_*\St$ is compatible with this decomposition, and its restriction to $\M$ coincides with $S$ (by using (A.3.1) below together with the remark after (1.2.2)).
We then get the direct sum decompositions
$$H^j(\Xt,\MMt)=H^j(X,\M)\oplus H^{j-2}(Z,\M_Z)(-1)\q\q(j\in\Z).
\leqno(2.5.2)$$
where $H^j(X,\M):=H^j(a_X)_*\M$ for the structure morphism $a_X:X\to pt$, and similarly for $H^j(\Xt,\MMt)$, etc.
The above argument implies that these direct sum decompositions are compatible with the induced pairing by $\St$, and moreover its restriction to the first factors coincides with the induced pairing by $S$.
\sk
We have the Lefschetz pencil
$$p:\Xt\to\PP^1,$$
and Theorem~(1.4) for $f(Z)\ne pt$ can be applied to this.
So the proof of Theorem~(1.4) for $\dim Z=d$ and $f(Z)=pt$ is reduced to the case $X=\PP^1$, where we can apply Zucker's result \cite{Zu}.
(Note that Zucker gave an ``algebraic description" of the Hodge filtration $F$ using holomorphic differential forms, see \cite[Corollary 6.15 and Proposition 9.1]{Zu}.)
However, we have to make some more calculations about polarizations on primitive classes related to the Leray spectral sequence, etc.\ (which are not quite trivial, see \cite[Sections 5.3.8\,-11]{mhp} for details).
\sk
The Lefschetz pencil is also used in an essential way for the proof of the Weil conjecture and the hard Lefschetz theorem in \cite{weil}.
The reduction to the curve case is by analogy with the $\ell$-adic case in some sense.
\msn
{\bf Remarks.} (i) For the calculation of the nearby and vanishing cycle functors in the normal crossing case, we use the so-called ``combinatorial description" of mixed Hodge modules with normal crossing singular loci. Here the {\it compatibility} of the $d_X+1$ filtrations $F,V_{(i)}$ ($i\in[1,d_X])$ is quite essential, where $V_{(i)}$ is the $V$-filtration along $x_i=0$, and the $x_i$ are local coordinates compatible with the singular locus of a mixed Hodge module. Note, however, that this description {\it never} gives an equivalence of categories (consider, for instance, the case of variations of mixed Hodge structure having no singular loci; in fact, this ``Hodge-theoretic combinatorial description" gives only the information of the fiber at the fixed point). Nevertheless it is quite useful when it is combined with the Verdier-type extension theorem \cite{Ve2} inductively, see also \cite[Proposition 3.13]{mhm}, etc.
\sk
It seems rather easy to predict Hodge-theoretic combinatorial formulas for the nearby and vanishing cycle functors together with the induced pairing in the normal crossing case. These are implicitly related with Beilinson's construction of nearby cycles \cite{Bei2}, see \cite{ext} for the mixed Hodge module case. It seems more difficult to prove that these formulas actually hold, see for instance \cite{dual}.
(Note that the argument in the Appendix of \cite{mhm} was simplified by the writer. The original argument used the reduction to the 2-dimensional case, and was much more complicated.)
\sk
(ii) The results of Cattani, Kaplan, Schmid (\cite{CK}, \cite{CKS1}, \cite{CKS2}, \cite{CKS3}) are used in an essential way for the above ``Hodge-theoretic combinatorial description". For instance, the {\it descent lemma} in \cite{CKS1}, \cite{CKS3} is crucial to the ``combinatorial description" of the pure Hodge module corresponding to the intersection complex. (This lemma is called ``the vanishing cycle theorem" in \cite{KK4}, which does not seem to be contained in \cite{KK2}.)
\sk
(iii) It is still an open problem whether the Hodge structure obtained by the $L^2$ method in \cite{CKS3}, \cite{KK4} has an ``algebraic description" using holomorphic differential forms, and in the algebraic case, whether it coincides with the Hodge structure obtained by the theory of mixed Hodge modules. (It has been expected that the detailed versions of \cite{KK3}, \cite{KK5} would give a positive answer to these problems, and some people thought that \cite{KK2}, \cite{KK4} were written for this purpose. As for \cite{BrSY}, it seems rather difficult to apply it to filtered $L^2$-sheaf complexes.) We have to assume that polarizable variations of Hodge structure are geometric ones in \cite{toh}, and also in \cite{PS} for the analytic case.
\sk
(iv) In the curve case, the answer to the above first problem was already given in \cite[Corollary 6.15 and Proposition 9.1]{Zu}, and the second problem is then easy to solve, see \cite[Section 5.3.10]{mhp}.
\sk
(v) It does not seem easy to generalize the results in \cite{CKS3}, \cite{KK4} to the case of a ``tubular neighborhood" of a subvariety in a smooth complex algebraic variety, since we would have to take a {\it complete} metric to get a Hodge structure by applying a standard method.
\msn
{\bf 2.6.~Outline of proof of Theorem~(2.4).}
We have to study the weight spectral sequence
$$E_1^{-k,j+k}=\Hc^jf_*\Gr_k^W\psi_h\M\Longrightarrow\Hc^jf_*\psi_h\M,
\leqno(2.6.1)$$
and similarly with $\psi_h$ replaced by $\varphi_{h,1}$.
This spectral sequence is defined in $\MF(Y,\Q)$ if the differentials $d_r$ ($r\ges 1$) are strictly compatible with the filtration $F$ inductively, see \cite[1.3.6]{mhp}. The last assertion for $r=1$ follows from (2.3.3), since $d_1$ preserves the Hodge filtration $F$ and also the weight filtration $W$ (the latter is shifted depending on the degree $j$ of the cohomology sheaf $\Hc^jf_*\psi_h\M$). For $r\ges 2$, we can show $d_r=0$ inductively by using (2.3.4). In particular, the $E_2$-degeneration of the spectral sequence follows. The above argument also implies the strictness of the filtration $F$ on the direct image $f_*^{\D}\psi_h\M$ by using the theory of compatible filtrations ({\it loc.~cit.}). A similar assertion holds for $\varphi_{g,1}\M$. These assertions imply the assertion~(i) of Theorem~(2.4) by using the completion by the $V$-filtration, see \cite[Section 3.3]{mhp} for details.
\sk
To show the remaining assertions, we have to show that the ``bi-symmetry" for the actions of $\ell$, $N$ on the $E_1$-term is preserved on the $E_2$-terms, and moreover the ``bi-primitive part" of the $E_2$-term is {\it represented} by the bi-primitive part of the $E_1$-term. Using the strict support decomposition together with the easy part of Theorem~(1.3), we can reduce the assertions to the case where the spectral sequence is defined in the category of Hodge structures. Theorem~(2.4) then follows from the theory of bi-graded Hodge structures of Lefschetz type as in the proof of \cite[Proposition 4.2.2]{mhp} (see also \cite{pos} where a slightly better explanation is given).
\msn
{\bf Remarks.} (i) Signs were not determined in \cite[Proposition 4.2.2]{mhp}, since this was a very subtle issue at that time (see for instance \cite[Section 2.2.5]{th2} where a problem of sign for Chern classes was raised.) In the case of the nearby cycles of the constant sheaf $\Q_X$ in the normal crossing case (that is, in the case of Steenbrink \cite{St}), the primitive part
$${}^P\Gr^W_{d_X-1+i}\psi_{h,1}(\Q_X[d_X-1])$$
is the direct sum of the constant sheaves supported on intersections of $i+1$ irreducible components of $h^{-1}(0)$. In particular, its support is pure dimensional, and the signs should depend only on this dimension. Then there would be no problem for using \cite[Proposition 4.2.2]{mhp} in this case. In the general case, however, the situation is much more complicated. In fact, there may be direct factors $({}^P\Gr^W_k\psi_h\M)_Z$ of the primitive part ${}^P\Gr^W_k\psi_h\M$ which have strict supports $Z$ of {\it various dimensions}, and surject to the {\it same} closed subvariety of $Y$ (for instance, if $g^{-1}(0)=\{0\}$). In this case, we have to determine exactly the sign for each direct factor so that the {\it positivity} becomes compatible among the direct images of direct factors with various strict support dimensions.
\sk
(ii) Precise signs are written in \cite[Lemma 5.3.6]{mhp} following Deligne's sign system \cite{sign} (see also \cite{GN}). Note that the conclusion of Lemma 5.3.6 follows from the {\it proof} of \cite[Proposition 4.2.2]{mhp}, since the hypothesis of the lemma is stronger than that of the proposition and the conclusion is essentially the same (that is, the $E_2$-term is bi-symmetric for $\ell,N$ and its bi-primitive part is represented by the bi-primitive part of the $E_1$-term). A pairing on the $E_1$-term of the Steenbrink weight spectral sequence \cite{St} satisfying Deligne's sign system \cite{sign} is constructed in \cite{GN} although its relation with the pairing induced by the nearby cycle functor does not seem to be clear (for instance, an isomorphism like (1.4.10) does not seem to be used there).
\sk
(iii) An argument showing the decomposition (2.1.11) is noted in \cite[Lemma 5.2.5]{mhp} which can replace \cite[Corollary 4.2.4]{mhp} if one is quite sure that the signs in the lemma really hold in the case one is considering. In fact, the assertion is very sensitive to the signs: if the signs are modified, then the role of $H$ and $H'$ can be reversed, and we may get a decomposition of $H$ instead of $H'$. (There is a misprint in the last line of the lemma: $S'$ should be $H'$.) If one is not completely sure whether the signs of the lemma really hold in the case under consideration, then it is still possible to use \cite[Corollary 4.2.4]{mhp} at least in the constant sheaf case with normal crossing singularities as is explained in Remark~(i) above.
\bs\bs
\vbox{\centerline{\bf 3. Mixed Hodge modules}
\bsn
In this section we explain a simplified definition of mixed Hodge modules following \cite{def}.}
\msn
{\bf 3.1.~Admissibility condition for weakly mixed Hodge modules.}
Let $X$ be a smooth complex algebraic variety, and $g\in\Gamma(X,\Oc_X)$. For a weakly mixed Hodge module
$$(\M,W)=((M;F,W),(K,W))\in\MHW(X),$$
(see (2.3.1)), set
$$(\Mt;F,W)=(i_g)_*^{\D}(M;F,W),$$
where $W$ is not shifted under the direct image. We have the filtration $V$ on $\Mt$ as in (B.6) below. We define the filtration $L$ on the nearby and vanishing cycle functors by
$$L_k\psi_g\M:=\psi_gW_{k+1}\M,\q L_k\varphi_{g,1}\M=\varphi_{g,1}W_k\M\q(\forall\,k\in\Z).
\leqno(3.1.1)$$
Here the filtration $F$ can be neglected when the filtration $L$ is defined.
\sk
We say that $\M$ is {\it admissible} along $g=0$ (or $g$-{\it admissible} for short) if the following two conditions are satisfied:
$$\h{Three filtrations $F,W,V$ on $\Mt$ are compatible filtrations (see (C.2) below).}
\leqno(3.1.2)$$
\vskip-7mm
$$\h{There is the relative monodromy filtration $W$ on $(\psi_g\M,L)$, $(\varphi_{g,1}\M,L)$.}
\leqno(3.1.3)$$
The last condition means that there is a unique filtration $W$ on $\psi_g\M$ satisfying the following two conditions:
$$\aligned N(W_i\psi_g\M)&\subset W_{i-2}\psi_g\M(-1)\q(\forall\,i\in\Z),\\ N^i:\Gr^W_{k+i}\Gr^L_k\psi_g\M&\simto\Gr^W_{k-i}\Gr^L_k\psi_g\M(-i)\q(\forall\,i\in\N,\,k\in\Z),\endaligned
\leqno(3.1.4)$$
and similarly for $\varphi_{g,1}\M$, see \cite{weil}, \cite{StZ}.
(Here the filtration $F$ can be forgotten. However, the last isomorphism is compatible with $F$ by (2.3.3) if $(\psi_f\M,W),(\varphi_{f,1}\M,W)\in\MHW(X)$.)
\msn
{\bf Remark.} Let $X$ be a smooth complex variety, and $\Xo$ be a smooth compactification such that $D:=\Xo\setminus X$ is a divisor with simple normal crossings. Let
$$(\M,W)=\bl((M;,F,W),(H,W)\br)$$
be a variation of mixed Hodge structure on $X$ such that $\Gr_k^W\M$ are polarizable pure Hodge structures of weight $k$ for any $k\in\Z$, where $(M;F,W)$ is the underlying bi-filtered $\Oc_X$-module, and $(H,W)$ is the underlying filtered $\Q$-local system.
\sk
Assume the local monodromies of $H$ are all {\it unipotent}. Let $\Mo$ be the {\it canonical}\, Deligne extension of $M$ over $\Xo$ (that is, the residues of the logarithmic connections are all {\it nilpotent}, see \cite{eq}). The filtrations $F,W$ are naturally extended on $\Mo$ by taking the intersection with $\Mo$ of the open direct images of $F,W$ under the inclusion $j:X\into\Xo$. Let $D_i$ be the irreducible components of $D$ ($i\in[1,r]$). Let $T_i$ be the local monodromy around a general point of $D_i$, which is defined on the fiber $H_{x_0}$ at a base point $x_0\in X$, and is unipotent by hypothesis. (This is well-defined up to a conjugate compatible with $W$.) Set $N_i:=(2\pi i)^{-1}\log T_i$. We denote by $L$ the filtration $W$ on $H_{x_0}$.
\sk
Under the above notation and assumption, $(\M,W)$ is an {\it admissible variation of mixed Hodge structure} in the sense of \cite{StZ}, \cite{Kadm} if and only if the following two conditions are satisfied:
$$\Gr_F^p\Gr_k^W\Mo\,\,\h{are locally free $\Oc_{\Xo}$-modules for any $p,k\in\Z$.}
\leqno(3.1.5)$$
\vskip-7mm
$$\h{There is the relative monodromy filtration on $(H_{x_0},L)$ for each $N_i$.}
\leqno(3.1.6)$$
(The last condition means that (3.1.4) holds with $\psi_g\M$ replaced by $H_{x_0}$.) In fact, this equivalence follows from \cite[Theorem 4.5.2]{Kadm}.
\sk
In the non-unipotent local monodromy case, let $\rho:X'\to X$ be a generically finite morphism of complex algebraic varieties such that $\rho^*\M$ has unipotent local monodromies around the divisor at infinity of a compactification of $X'$ (by replacing $X$ with a non-empty open subvariety if necessary). Then $\M$ is an admissible variation if and only if $\rho^*\M$ is. (In fact, may assume that $\rho$ is finite \'etale by shrinking $X$. Then we can take the direct image.)
\msn
{\bf 3.2.~Well-definedness of open direct images.}
Let $D$ be a locally principal divisor on a smooth complex algebraic variety $X$. Set $X':=X\setminus D$ with $j:X'\into X$ the inclusion.
We say that the open direct images $j_!,j_*$ are well-defined for $\M'\in\MHW(X')$, if there are
$$\aligned&\q\q\M'_!,\,\M'_*\in\MHW(X),\,\,\,\,\h{satisfying}\\{\rm rat}(\M'_!)&=j_!K',\,\,\,{\rm rat}(\M'_*)=\R j_*K'\q\h{with}\q K':={\rm rat}(\M'),\endaligned
\leqno(3.2.1)$$
(see (1.1.6) for rat), and moreover the following condition is satisfied:
$$\h{$\M'_!$, $\M'_*$ are $g$-admissible for any $g\in\Gamma(U,\Oc_U)$ with $g^{-1}(0)_{\rm red}=D_{\rm red}\cap U$.}
\leqno({\rm A3})$$
Here $U$ is any open subvariety of $X$.
If the above condition is satisfied, we then define
$$j_!\M':=\M'_!,\q j_*\M':=\M'_*.
\leqno(3.2.2)$$
If $\M'=j^{-1}\M$ with $\M\in\MHW(X)$ and $\M$ is $g$-admissible, then we have the canonical morphisms (see \cite[Proposition~2.11]{mhm})
$$j_!j^{-1}\M\to\M,\q\M\to j_*j^{-1}\M.
\leqno(3.2.3)$$
\msn
{\bf 3.3. Definition of mixed Hodge modules.} Let $X$ be a smooth complex algebraic variety. The category of mixed Hodge modules $\MHM(X)$ is the abelian full subcategory of $\MHW(X)$ in (2.3) defined by increasing induction on the dimension $d$ of the support as follows:
\sk
For $\M\in\MHW(X)$ with ${\rm supp}\,\M=Z$, we have $\M\in\MHM(X)$ if and only if, for any $x\in X$, there is a Zariski-open neighborhood $U_x$ of $x$ in $X$ together with $g_x\in\Gamma(U_x,\Oc_{U_x})$ such that
$$\dim Z\cap U_x\cap g_x^{-1}(0)<\dim Z,$$
$Z'_x:=Z\cap U_x\setminus g_x^{-1}(0)$ is smooth, and moreover the following two conditions are satisfied:
$$\h{$\M|_{Z'_x}$ is an admissible variation of mixed Hodge structure.}
\leqno(3.3.1)$$
$$\h{$\M|_{U_x}$ is $g_x$-admissible, and $\varphi_{g_x,1}\M|_{U_x}\in\MHM(U_x)$.}
\leqno(3.3.2)$$
More precisely, (3.3.1) means that $\M|_{U'_x}$ is isomorphic to the direct image of an admissible variation of mixed Hodge structure on $Z'_x$ by the closed embedding
$$i_{Z'_x}:Z'_z\into U'_x:=U_x\setminus g_x^{-1}(0).$$
\ms
If $Z=\{x\}$ for $x\in X$, then we set
$$\MHM_{\{x\}}(X):=\MHW_{\{x\}}(X)={\rm MHS}(\Q),
\leqno(3.3.3)$$
where the first and second categories are respectively full subcategories of $\MHM(X)$ and $\MHW(X)$ consisting of objects supported on $x$, and the last one is the category of graded-polarizable mixed $\Q$-Hodge structures \cite{th2}. (Here the direct image by $\{x\}\into X$ is used.)
\ms
This definition is justified by the following (see \cite[Theorem 1]{def}).
\msn
{\bf Theorem~3.4.} {\it Conditions {\rm (3.3.1--2)} are independent of the choice of $U_x$, $g_x$. More precisely, if they are satisfied for some $U_x$, $g_x$ for each $x\in Z$, then $(3.3.2)$ is satisfied for any $U_x$, $g_x$, and $(3.3.1)$ is satisfied in case ${\rm rat}(\M')$ is a local system up to a shift of complex, see $(1.1.6)$ for ${\rm rat}$.}
\ms
We have moreover the following (see \cite[Theorem 2]{def}).
\msn
{\bf Theorem~3.5.} {\it The categories $\MHM(X)$ for smooth complex algebraic varieties $X$ are stable by the canonically defined cohomological functors $\Hc^jf_*$, $\Hc^jf_!$, $\Hc^jf^*$, $\Hc^jf^!$, $\psi_g$, $\varphi_{g,1}$, $\boxtimes$, $\DD$, where $f$ is a morphism of smooth complex algebraic varieties and $g\in\Gamma(X,\Oc_X)$. Moreover these functors are compatible with the corresponding functors of the underlying $\Q$-complexes via the forgetful functor {\rm rat} in $(1.1.6)$.}
\ms
The proofs of these theorems use Beilinson's maximal extension together with the stability by subquotients systematically. The well-definedness of open direct images in (3.3) is reduced to the normal crossing case, see \cite{def}.
Combining Theorem~(3.5) with the construction in \cite{mhm}, we can get the following (see \cite[Corollary 1]{def}).
\msn
{\bf Theorem~3.6.} {\it There are canonically defined functors $f_*$, $f_!$, $f^*$, $f^!$, $\psi_g$, $\varphi_{g,1}$, $\boxtimes$, $\DD$, $\otimes$, ${\Hc}om$ between the bounded derived categories $D^b\MHM(X)$ for smooth complex algebraic varieties $X$ so that we have the canonical isomorphisms $H^jf_*=\Hc^jf_*$, etc., where $f$ is a morphism of smooth complex algebraic varieties, $g\in\Gamma(X,\Oc_X)$, $H^j$ is the usual cohomology functor of the derived categories, and $\Hc^jf_*$, etc.\ are as in Theorem~$(3.5)$. Moreover the above functors between the $D^b\MHM(X)$ are compatible with the corresponding functors of the underlying $\Q$-complexes via the forgetful functor {\rm rat}.}
\msn
{\bf 3.7. Some notes on references about applications.} Since there is no more space to explain about applications of mixed Hodge modules, we indicate some references here. These are not intended to be complete.
\sk
For applications to algebraic cycles, see
\cite{BRS}, \cite{BFNP}, \cite{MuS}, \cite{NS}, \cite{RS1}, \cite{RS2}, \cite{int}, \cite{HC1}, \cite{HC2}, \cite{HTC}, \cite{ari}, \cite{rcy}, \cite{FCh}, \cite{Dir}, \cite{Tho}, \cite{HCh}, \cite{SS2}, etc.
Some of them are related to normal functions. For the latter, see also
\cite{ext}, \cite{anf}, \cite{spr}, \cite{SS1}, \cite{Schn}, etc.
\sk
Related to mixed Hodge structures on cohomologies of algebraic varieties, see
\cite{BDS}, \cite{DS1}, \cite{DS2}, \cite{DS3}, \cite{DS5}, \cite{DS6}, \cite{DSW}, \cite{DuS}, \cite{OS}, \cite{PS}, \cite{cmp}, \cite{Fil}, \cite{SZ}, etc.
\sk
About Bernstein-Sato polynomials, Steenbrink spectra, and multiplier ideals, see
\cite{Bu}, \cite{BMS}, \cite{BS1}, \cite{BS2}, \cite{BSY}, \cite{DMS}, \cite{DMST}, \cite{DS4}, \cite{DS7}, \cite{MP}, \cite{ste}, \cite{rat}, \cite{mic}, \cite{bhyp}, \cite{Mul}, \cite{pow}, etc.
\sk
For direct images of dualizing sheaves and vanishing theorems, see \cite{FFT}, \cite{kol}, \cite{Su}, etc.
Concerning Hirzebruch characteristic classes, see 
\cite{BrScY}, \cite{MS}, \cite{MSS1}, \cite{MSS2}, \cite{MSS3}, \cite{MSS4}, etc.
\bs\bs
\vbox{\centerline{\bf Appendix A. Hypersheaves}
\bsn
In this appendix we review some basics of hypersheaves, see Convention~3.}
\msn
{\bf A.1.} Let $X$ be a complex algebraic variety or a complex analytic space, and $A$ be a subfield of $\C$. We denote by $D_c^b(X,A)$ the derived category of bounded complexes of $A$-modules with constructible cohomology sheaves, see \cite{Ve1}, etc.
In the algebraic case, we use the classical topology for the sheaf complexes although we assume that {\it stratifications are algebraic}.
\sk
The category of hypersheaves $\HS(X,A)$ is the {\it full subcategory} of $D_c^b(X,A)$ consisting of objects $K$ satisfying the condition:
$$\dim{\rm supp}\,\Hc^iK\les-i,\q\dim{\rm supp}\,\Hc^i\1\DD(K)\les-i\q(\forall\,i\in\Z).
\leqno({\rm A}.1.1)$$
Here $\Hc^iK$ is the $i$\,th cohomology sheaf of $K$ in the usual sense, and $\DD(K)$ is the dual of $K$. The latter can be defined by
$$\DD(K):=\R\Hc om_A(K,\DD_X),
\leqno({\rm A}.1.2)$$
with $\DD_X$ the dualizing sheaf in $D^b_c(X,A)$. In the {\it smooth} case (by taking an embedding into smooth varieties), it can be defined by
$$\DD_X:=A_X(d_X)[2d_X],
\leqno({\rm A}.1.3)$$
with $d_X:=\dim X$.
\sk
By \cite{BBD}, $\HS(X,A)$ is an {\it abelian category}, and there are canonical cohomological functors
$$\mh\Hc^i:D^b_c(X,A)\to\HS(X,A)\q(i\in\Z),
\leqno({\rm A}.1.4)$$
where the superscript $\,\mh\,$ means the ``middle perversity".
\msn
{\bf A.2.~Nearby and vanishing cycles.} Let $g$ be a holomorphic function on an analytic space $X$. Let $\Delta\subset\C$ be a sufficiently small open disk with center $0$, and $\pi:\widetilde{\Delta^*}\to\Delta^*$ be a universal covering of the punctured disk $\Delta^*$. Let $\pi':\widetilde{\Delta^*}\to\C$ be its composition with the inclusion $\Delta^*\into\C$.
Let $X_{\infty}$ be the base change of $X$ by $\pi'$. We denote by $\widetilde{j}:X_{\infty}\to X$ the base change of $\pi'$ by $g$. Set $X_0:=f^{-1}(0)$ with $i_0:X_0\into X$ the canonical inclusion. The nearby and vanishing cycle functors
$$\psi_g,\,\varphi_g:D_c^b(X,A)\to D_c^b(X_0,A)$$
are defined as in \cite{van} by
$$\psi_gK:=i_0^*\,\R\widetilde{j}_*\widetilde{j}^*K,\,\,\,\varphi_gK:=C(i_0^*\,K\to\psi_gK)\q\h{for}\,\,\,\,K\in D_c^b(X,A),
\leqno({\rm A}.2.1)$$
where we take a flasque resolution of $K$ to define $\R\widetilde{j}_*\widetilde{j}^*K$ and also the mapping cone. By definition we have a distinguished triangle
$$i^*K\to\psi_gK\to\varphi_gK\buildrel{+1}\over\to.
\leqno({\rm A}.2.2)$$
\sk
The action of the monodromy $T$ is defined by $\gamma^*$ with $\gamma$ the automorphism of $\widetilde{\Delta^*}$ over $\Delta^*$ defined by $z\mapsto z+1$. Here $\widetilde{\Delta^*}$ is identified with $\{z\in\C\mid{\rm Im}\,z>r\}$ for some $r>0$ and $\pi$ is given by $z\mapsto t:=\exp(2\pi iz)$. (This is compatible with the usual definition of the monodromy of a local system $L$ on $\Delta^*$. In fact, $(\gamma^*\sigma)(z_0)=\sigma(z_0+1)$ for $\sigma\in\Gamma(\widetilde{\Delta^*},\pi^*L)$ with $z_0$ a base point of $\widetilde{\Delta^*}$, and the monodromy is given by the composition of canonical isomorphisms: $L_{\pi(z_0)}=(\pi^*L)_{z_0}=\Gamma(\widetilde{\Delta^*},\pi^*L)=(\pi^*L)_{z_0+1}=L_{\pi(z_0)}$.)
There is a nonzero minimal polynomial for $T$ locally on $X$, and this implies the Jordan decomposition $T=T_sT_u$ (with $T_s,T_u$ polynomials in $T$ locally on $X$).
\sk
Assume $K\in\HS(X,A)$. Set
$$\mh\psi_gK:=\psi_gK[-1],\,\,\mh\varphi_gK:=\varphi_gK[-1].
\leqno({\rm A}.2.3)$$
Then
$$\mh\psi_gK,\,\mh\varphi_gK\in\HS(X_0,A).
\leqno({\rm A}.2.4)$$
This follows for instance from \cite{Kvan}, \cite{Ma3} by using the Riemann-Hilbert correspondence.
\sk
In the case $A=\C$, this implies the decompositions in the abelian category $\HS(X_0,A)$:
$$\mh\psi_gK=\mopl_{\la\in\C^*}\,\mh\psi_{g,\la}K,\q\mh\varphi_gK=\mopl_{\la\in\C^*}\,\mh\varphi_{g,\la}K,
\leqno({\rm A}.2.5)$$
(which are locally finite direct sum decompositions), where
$$\mh\psi_{g,\la}K:={\rm Ker}(T_s-\la)\subset\mh\psi_gK,\q\mh\varphi_{g,\la}K:={\rm Ker}(T_s-\la)\subset\mh\varphi_gK.
\leqno({\rm A}.2.6)$$
\sk
In the case $A\subset\C$, we have only the decompositions
$$\mh\psi_gK=\mh\psi_{g,1}K\oplus\mh\psi_{g,\ne 1}K,\q\mh\varphi_gK=\mh\varphi_{g,1}K\oplus\mh\varphi_{g,\ne 1}K,
\leqno({\rm A}.2.7)$$
which are compatible with the above decompositions after taking the scalar extension by $A\into\C$.
\sk
By (A.2.2) we have the canonical isomorphisms
$$\mh\psi_{g,\ne1}K\simto\mh\varphi_{g,\ne1}K,\q\mh\psi_{g,\la}K\simto\mh\varphi_{g,\la}K\,\,\,\,(\la\ne 1,\,A=\C),
\leqno({\rm A}.2.8)$$
since the action of $T$ on $i^*K$ is trivial.
\sk
If $K=A_X$ and $X$ is a smooth algebraic variety (or a complex manifold with $X_0$ compact), then the nearby cycle functor $\psi_gA_X$ is also defined by
$$\psi_gA_X=\R\rho_*A_{X_c},
\leqno({\rm A}.2.9)$$
where $X_c:=f^{-1}(c)\subset X$ with $c\in\C^*$ sufficiently near $0$, and $\rho:X_c\to X_0$ is an appropriate contraction morphism. The latter is constructed by using an embedded resolution of $X_0\subset X$.
\msn
{\bf A.3.~Compatibility with the dual functor $\DD$.} We say that a pairing
$$K\otimes_AK'\to\DD_X(k)$$
is a {\it perfect pairing} (with $k\in\Z$) if its corresponding morphism
$$K\to\DD(K')(k)=\R\Hc om_A(K',\DD_X)(k)$$
is an isomorphism in $D^b_c(X,A)$. Here $\DD_X$ is as in (A.1.2), and the above correspondence comes from the isomorphism
$${\rm Hom}(K\otimes K',\DD_X(k))={\rm Hom}(K,\R\Hc om_A(K',\DD_X(k))).
\leqno({\rm A}.3.1)$$
The Tate twist $(k)$ is defined as in \cite[Definition 2.1.13]{th1}, see also a remark after (2.1.10). This can be neglected if $i=\sqrt{-1}$ is chosen. However, this twist is quite useful in order to keep track of ``weight". In fact, if $K=K'$ and it has pure weight $w$, then $\DD K$ should have weight $-w$, and the above $k$ must be equal to $-w$, since the Tate twist $(k)$ changes the weight by $-2k$ ({\it loc.~cit.}).
\sk
Assume there is a perfect pairing
$$S:K\otimes_A K'\to\DD_X(-w)\q\h{for}\,\,\,K,K'\in\HS(X,A).$$
It induces a canonical pairing
$$\psi_gS:\psi_gK\otimes_A\psi_gK'\to\psi_g\DD_X(-w).
\leqno({\rm A}.3.2)$$
\sk
Assume $X$ is {\it smooth} (by taking an embedding into a smooth variety), and $X_0$ is also {\it smooth} (by replacing $K$ with its direct image under the graph embedding by $g$). Then we have
$$\psi_g\DD_X=A_{X_0}(d_X-w)[2d_X]=\DD_{X_0}(1-w)[2],
\leqno({\rm A}.3.3)$$
and (A.3.2) indues a canonical perfect pairing
$$\mh\psi_gS:\mh\psi_gK\otimes_A\mh\psi_gK'\to\DD_{X_0}(1-w).
\leqno({\rm A}.3.4)$$
Here some sign appears, and this is closely related to the sign in (1.4.8).
\sk
The above construction is compatible with the monodromy $T$, that is,
$$\mh\psi_gS=\mh\psi_gS\ssc(T\otimes T).$$
Since $T^e$ is unipotent for some $e\in\Z_{>0}$, this implies
$$\mh\psi_gS\ssc(N\otimes id)=-\mh\psi_gS\ssc(id\otimes N),\q\mh\psi_gS=\mh\psi_gS\ssc(T_s\otimes T_s),
\leqno({\rm A}.3.5)$$
where $T=T_sT_u$ is the Jordan decomposition, and $N:=(2\pi i)^{-1}\log T_u$.
\sk
We then get the induced perfect pairings
$$\aligned\mh\psi_{g,1}S:\mh\psi_{g,1}K&\otimes_A\mh\psi_{g,1}K'\to\DD_{X_0}(1-w),\\\mh\psi_{g,\ne1}S:\mh\psi_{g,\ne1}K&\otimes_A\mh\psi_{g,\ne1}K'\to\DD_{X_0}(1-w),\\\mh\psi_{g,\la}S:\mh\psi_{g,\la}K&\otimes_A\mh\psi_{g,\la^{-1}}K'\to\DD_{X_0}(1-w)\q(A=\C).\endaligned
\leqno({\rm A}.3.6)$$
\sk
For the vanishing cycle functor $\varphi_g$, we have the induced perfect pairing
$$\mh\varphi_{g,1}S:\mh\varphi_{g,1}K\otimes_A\mh\varphi_{g,1}K'\to\DD_{X_0}(-w),
\leqno({\rm A}.3.7)$$
satisfying
$$\mh\varphi_{g,1}S\ssc({\rm can}\otimes id)=\mh\psi_{g,1}S\ssc(id\otimes{\rm Var}),
\leqno({\rm A}.3.8)$$
where the morphisms
$${\rm can}:\mh\psi_{g,1}K\to\mh\varphi_{t,1}K,\q{\rm Var}:\mh\varphi_{g,1}K'\to\mh\psi_{t,1}K'(-1),
\leqno({\rm A}.3.9)$$
are constructed in \cite[Section 5.2.1]{mhp}. (These correspond respectively to the morphisms $-\Gr_V\dd_t$, $\Gr_Vt$ in (B.6.9) below if $K=\DR_X(M)$, $K'=\DR_X(M')$ with $A=\C$.) We have
$${\rm Var}\ssc{\rm can}=N,\,\,\,{\rm can}\ssc{\rm Var}=N.
\leqno({\rm A}.3.10)$$
\sk
Note that the target of (A.3.7) is different from that of (A.3.6) by the Tate twist, and $\mh\varphi_{g,\ne 1}S$ is given by $\mh\psi_{g,\ne 1}S$ together with the isomorphism (A.2.8).
The construction of (A.3.7) is not quite trivial, see \cite[Sections 5.2.1 and 5.2.3 and Lemma 5.2.4]{mhp}. For instance, we used there the isomorphism
$$i^!K'=[i^*K'\to\psi_{g,1}K'\buildrel{-N}\over\longrightarrow\psi_{g,1}K'(-1)],
\leqno({\rm A}.3.11)$$
where $i:X_0\into X$ is the inclusion. If we replace this with
$$i^!K'=[i^*K'\to\psi_{g,1}K'\buildrel{id-T}\over\longrightarrow\psi_{g,1}K'],
\leqno({\rm A}.3.12)$$
then (A.3.8) would hold with Var, $N$ respectively replaced by var, $T-id$, where the Tate twist should be omitted (since $T-id$ is not compatible with the weight structure).
\bs\bs
\vbox{\centerline{\bf Appendix B. $\D$-modules}
\bsn
In this appendix we review some basics of $\D$-modules.}
\msn
{\bf B.1.~Holonomic $\D$-modules.} Let $X$ be a complex manifold of dimension $d_X$, and $M$ be a coherent left $\D_X$-module. This means that $M$ has locally a finite presentation
$$\mopl^p\,\D_U\to\mopl^q\,\D_U\to M|_U\to 0,$$
over sufficiently small open subsets $U\subset X$.
(This is equivalent to the condition that $M$ is quasi-coherent over $\Oc_X$ and is locally finitely generated over $\D_X$.)
\sk
A filtration $F$ on $M$ is called a {\it good filtration} if $(M,F)$ satisfies
$$(F_p\D_X)\1(F_qM)\subset F_{p+q}M\q(p\in\N,\,q\in\Z),
\leqno({\rm B}.1.1)$$
and $\Gr^F_{\ssb}M$ is a {\it coherent} $\Gr^F_{\ssb}\D_X$-module.
(The last condition is equivalent to the conditions that each $F_pM$ is coherent over $\Oc_X$ and the equality holds for $q\gg0$ in (B.1.1).)
Here $F$ on $\D_X$ is by the order of differential operators, that is, we have for local coordinates $(x_1,\dots,x_{d_X})$
$$F_p\D_X=\msum_{|\nu|\les p}\Oc_X\dd_{x_i}^{\nu_i}.$$
\sk
The {\it characteristic variety} ${\rm Ch}(M)\subset T^*X$ of a coherent left $\D_X$-module $M$ is defined to be the support of the $\Gr^F_{\ssb}\D_X$-module $\Gr^F_{\ssb}M$ in the cotangent bundle $T^*X$. (The latter can be defined to be the union of the analytic subspaces of $T^*X$ defined by the ideal of $\Gr^F_{\ssb}\D_X$ annihilating $g_i$ with $g_i$ local generators of $\Gr^F_{\ssb}M$. Here $\Gr^F_{\ssb}\D_X$ is identified with the sheaf of holomorphic functions on $T^*X$ which are polynomials on fibers of $T^*X\to X$.) This is independent of a choice of a good filtration $F$.
\sk
By the involutivity of the characteristic varieties (see \cite{SKK}, \cite{Ma2}, \cite{Ga}), it is known that
$$\dim{\rm Ch}(M)\ges\dim X.
\leqno({\rm B}.1.2)$$
(See also \cite{Bo} for the algebraic $\D$-module case.)
\sk
A coherent left $\D_X$-module $M$ is called {\it holonomic} if
$$\dim{\rm Ch}(M)=\dim X.
\leqno({\rm B}.1.3)$$
We will denote by $M_{\rm hol}(\D_X)$ the abelian category of holonomic $\D_X$-modules.
\msn
{\bf B.2.~Left and right $\D$-modules.}
We have the transformation between filtered left and right $\D_X$-modules on a complex manifold $X$ of dimension $d_x$ which associates the following to a filtered left $\D_X$-module $(M,F)$:
$$(\Om_X^{d_X},F)\otimes_{\Oc_X}(M,F),
\leqno({\rm B}.2.1)$$
where the filtration $F$ on $\Om_X^{d_X}$ is defined by the condition
$$\Gr^F_p\Om_X^{d_X}=0\q(p\ne-d_X).
\leqno({\rm B}.2.2)$$
So the filtration is shifted by $-d_X$.
Here it is better to distinguish $\Om_X^{d_X}$ and the dualizing sheaf $\omega_X$, since the Hodge filtration $F$ on $\omega_X$ is usually defined by
$$\Gr^F_p\omega_X=0\q(p\ne 0).
\leqno({\rm B}.2.3)$$
\sk
By choosing local coordinates $x_1,\dots,x_{d_X}$, the sheaf $\Om_X^{d_X}$ is trivialized by $\ddd x_1\wedge\cdots\wedge \ddd x_{d_X}$ locally on $X$, and forgetting $F$, the transformation is given by the anti-involution $^*$ of $\D_X$ defined by the conditions (see for instance \cite{Ma1}):
$$\dd_{x_i}^*=-\dd_{x_i},\q g^*=g\,\,\,(g\in\Oc_X),\q(PQ)^*=Q^*P^*\,\,\,(P,Q\in\D_X).
\leqno({\rm B}.2.4)$$
\sk
For a right $\D_X$-module $N$, the left $\D_X$-module corresponding to it is denoted often by
$$N\otimes_{\Oc_X}(\Om_X^{d_X})^{\vee}.
\leqno({\rm B}.2.5)$$
Here $L^{\vee}$ denotes the dual of a locally free sheaf $L$ in general, that is, $L^{\vee}:=\Hc om_{\Oc_X}(L,\Oc_X)$.
\msn
{\bf B.3.~Direct images.} For a closed embedding $i:X\to Y$ of smooth complex algebraic varieties, the direct image of a filtered {\it right} $\D$-module $(M,F)$ is defined by
$$i_*^{\D}(M,F)=(M,F)\otimes_{\D_X}(\D_{X\into Y},F),
\leqno({\rm B}.3.1)$$
where the sheaf-theoretic direct image is omitted to simplify the notation, and
$$(\D_{X\into Y},F):=\Oc_X\otimes_{\Oc_Y}(\D_Y,F).
\leqno({\rm B}.3.2)$$
\sk
For a filtered {\it left} $\D$-module $(M,F)$, the $\D$-module is twisted by $\omega_{X/Y}$ and the filtration $F$ is shifted by $r:~={\rm codim}_XY$ because of the transformation between filtered left and right $\D$-modules in (B.2).
If $X$ is locally defined by $y_1=\cdots=y_r$ with $y_1,\dots,y_m$ local coordinates of $Y$, then, setting $\dd_{y_i}:=\dd/\dd y_i$, the direct image is {\it locally} defined by
$$i_*^{\D}(M,F)=(M,F[r])\otimes_{\C}(\C[\dd_{y_1},\dots,\dd_{y_r}],F).
\leqno({\rm B}.3.3)$$
\sk
For a smooth projection $p:Z:=X\times Y\to Y$ with $X,Y$ smooth, the direct image of a filtered {\it left} $\D_Z$-module $(M,F)$ is defined by the sheaf-theoretic direct image of the relative de Rham complex $\DR_{Z/Y}(M,F)$, that is,
$$p_*^{\D}(M,F):=\R\1p_*\DR_{Z/Y}(M,F),
\leqno({\rm B}.3.4)$$
where $\DR_{Z/Y}(M,F)$ is the filtered complex defined by
$$(M,F)\to\Om_{Z/Y}^1\otimes_{\Oc_Z}(M,F[-1])\to\cdots\to\Om_{Z/Y}^{d_X}\otimes_{\Oc_Z}(M,F[-d_X]),
\leqno({\rm B}.3.5)$$
with the last term put at degree $0$. The differential of this complex is defined as in the absolute case (see (1.2.3)), and is locally given as the Koszul complex associated with the action of $\dd/\dd x_i$ on $M$ if $x_1,\dots,x_{d_X}$ are local coordinates of $X$.
(Note, however, that this does not work for {\it smooth morphisms} which are not necessarily {\it smooth projections}, since there is no canonical lift of vector fields on $Y$ to $Z$.)
\sk
In general, the direct image of a filtered right $\D_X$-module $(M,F)$ by a morphism of smooth varieties $f:X\to Y$ is defined by
$$i_*^{\D}(M,F)=p_*^{\D}\ssc(i_f)_*^{\D}(M,F),
\leqno({\rm B}.3.6)$$
where $i_f:X\to X\times Y$ is the graph embedding, and $p:X\times Y\to Y$ is the second projection so that $f=p\ssc i_f$.
\msn
{\bf Remarks.} (i) We can verify that the direct image in (B.3.6) is naturally isomorphic to the complex of the induced $\D_Y$-module associated with the (sheaf-theoretic) direct image of the filtered differential complex $\DR_X(M,F)$. This means the compatibility between the direct images of filtered differential complexes and filtered $\D$-modules.
\sk
(ii) It seems simpler to use the above construction of direct images instead of the induced $\D$-module construction as in \cite[3.3.6]{mhp} for the definition of direct images of $V$-filtrations.
\sk
(iii) The direct image for a morphism of singular varieties is rather complicated.
For the direct image of mixed Hodge modules, we may assume that the morphism is projective by using a Beilinson-type resolution (see \cite[Section 3]{Bei1} and the proof of \cite[Theorem 4.3]{mhm}), and the cohomological direct image is actually enough. So it is reduced to the case of a morphism of smooth varieties.
\msn
{\bf B.4.~Dual functor.} For a holonomic $\D_X$-module $M$ on a complex manifold $X$ of dimension $n$, its dual $\DD(M)$ is defined so that
$$(\Om_X^{d_X})^{\vee}\otimes_{\Oc_X}\DD(M)={\E}xt_{\D_X}^{d_X}(M,\D_X).
\leqno({\rm B}.4.1)$$
and $\DD$ is called the dual functor. This is a contravariant functor. It is well-known that
$$\E xt_{\D_X}^p(M,\D_X)=0\q(p\ne\dim X),
\leqno({\rm B}.4.2)$$
and
$$\DD^2=id.
\leqno({\rm B}.4.3)$$
\sk
We say that a filtered holonomic $\D_X$-module $(M,F)$ is {\it Cohen-Macaulay} if $\Gr_{\ssb}^FM$ is a Cohen-Macaulay $\Gr_{\ssb}^F$-module.
In this case, we have
$${\E xt}^i_{\Gr_{\ssb}^F\D_X}\bl(\Gr_{\ssb}^FM,\Gr_{\ssb}^F\D_X\br)=0\q\q(i\ne d_X),
\leqno({\rm B}.4.4)$$
and the dual filtered $\D_X$-module $\DD(M,F)$ can be defined so that
$$(\Om_X^{d_X},F)\otimes_{\Oc_X}\DD(M,F)=\R{\Hc om}_{\D_X}\bl((M,F),(\D_X,F[d_X])\br)[d_X].
\leqno({\rm B}.4.5)$$
This means that the last filtered complex is filtered quasi-isomorphic to a filtered $\D$-module.
\sk
It is known that the dual functor $\DD$ commutes with the de Rham functor $\DR_X$, that is, for a regular holonomic $\D_X$-module $M$, there is a canonical isomorphism (see for instance \cite[Proposition 2.4.12]{mhp}):
$$\DD(\DR_X(M))=\DR_X(\DD(M)).
\leqno({\rm B}.4.6)$$
\msn
{\bf B.5.~Regular holonomic $\D$-modules.} Let $Z$ be a closed analytic subset $Z$ of a complex manifold $X$. Let $\I_Z\subset\Oc_X$ be the ideal of $Z$. For a bounded complex of $\D_X$-modules $M^{\ssb}$, set
$$\aligned\Hc_{[Z]}^iM^{\ssb}&:=\rlap{\raise-10pt\h{$\,\,\scriptstyle k$}}\rlap{\raise-6pt\h{$\,\rightarrow$}}{\rm lim}\,\E{xt}^i_{\Oc_X}(\Oc_X/\I_Z^k,M^{\ssb})\,\bl(=\rlap{\raise-10pt\h{$\,\,\scriptstyle k$}}\rlap{\raise-6pt\h{$\,\rightarrow$}}{\rm lim}\,\E{xt}^i_{\D_X}(\D_X/\D_X\I_Z^k,M^{\ssb})\br),\\
\Hc_{[X|Z]}^iM^{\ssb}&:=\rlap{\raise-10pt\h{$\,\,\scriptstyle k$}}\rlap{\raise-6pt\h{$\,\rightarrow$}}{\rm lim}\,\E{xt}^i_{\Oc_X}(\I_Z^k,M^{\ssb})\,\bl(=\rlap{\raise-10pt\h{$\,\,\scriptstyle k$}}\rlap{\raise-6pt\h{$\,\rightarrow$}}{\rm lim}\,\E{xt}^i_{\D_X}(\D_X\I_Z^k,M^{\ssb})\br),\endaligned$$
so that we have a long exact sequence of $\D_X$-modules
$$\to\Hc_{[Z]}^iM^{\ssb}\to\Hc^iM^{\ssb}\to\Hc_{[X|Z]}^iM^{\ssb}\to\Hc_{[Z]}^{i+1}M^{\ssb}\to
\leqno({\rm B}.5.1)$$
see \cite{Khol2} (and also \cite{Gr} for the algebraic case).
Note that $\Hc_{[Z]}^0M$ for a holonomic $\D_X$-module $M$ is the largest holonomic $\D_X$-submodule supported in $Z$.
\ms
It is known that a holonomic $\D_X$-module $M$ with support $Z$ is {\it regular holonomic} if and only if there is a closed analytic subset $Z'\subset Z$ together with a proper morphism $\pi:\Zt\to Z$ such that $\dim Z'<\dim Z$, $Z\setminus Z'$ is smooth and equi-dimensional, $\Zt$ is smooth, $\pi$ induces an isomorphism over $Z\setminus Z'$, $\pi^{-1}(Z')$ is a divisor with normal crossings on $\Zt$, and moreover, by setting $\pi_X:=j_Z\ssc\pi$ with $j_Z:Z\into X$ the canonical inclusion, the following two conditions are satisfied:
\msn
(i) $\,\Hc^0_{[Z']}M$ is a regular holonomic $\D_X$-module.
\skn
(ii) $\,\Hc^0_{[X|Z']}M=\Hc^0(\pi_X)_*^{\D}\Mt$ with $\Mt$ the Deligne meromorphic extension of $M|_{Z\setminus Z'}$ on $\Zt$.
\ms
We may assume that $Z'$ is the union of ${\rm Sing}\,Z$ and the lower dimensional irreducible components of $Z$, and $\pi$ is given by the desingularization of the union of the maximal dimensional irreducible components of $Z$.
Note that in the case of algebraic $\D$-modules, we have $\Hc^0_{[X|Z']}M=j_*j^*M$ with $j:X\setminus Z'\into X$ the canonical inclusion.
\sk
The above criterion is by induction on the dimension of the support by using \cite{eq}, \cite{Khol2}. In fact, it is known that regular holonomic $\D$-modules are stable by subquotients and extensions in the category of holonomic $\D$-modules and also by the direct images under proper morphisms, and contain Deligne meromorphic extensions in the normal crossing case.
\ms
Let $M_{rh}(\D_X)\subset M_{\rm hol}(\D_X)$ be the full category of regular holonomic $\D_X$-modules on a complex manifold $X$. Let $D^b_{rh}(\D_X)\subset D^b(\D_X)$ be the full subcategory of bounded complexes with regular holonomic cohomologies. We have the equivalence of categories, that is, the Riemann-Hilbert correspondence:
$$\DR_X:D^b_{rh}(\D_X)\simto D_c^b(X,\C)\,\,\,\bl(\h{inducing}\,\,\,\DR_X:M_{rh}(\D_X)\simto\HS(X,\C)\br),
\leqno({\rm B}.5.2)$$
see \cite{Krh1}, \cite{Krh2}, \cite{KK1}, \cite{Mthe}, \cite{Mrh}, \cite{Meq} (and also \cite{Bo}, \cite{HT}, etc.\ for the algebraic case).
\msn
{\bf Remarks.} (i) There are many ways to define the full subcategory $M_{rh}(\D_X)\subset M_{\rm hol}(\D_X)$. Their equivalences may follow by using the above argument in certain cases.
\sk
(ii) There is a nontrivial point about the commutativity of some diagram used in certain proofs of (B.5.2), and this is studied in \cite[Section 4]{ind}.
\sk
(iii) Some people say that (B.5.2) is essentially proved in \cite{KK1} where the full faithfulness of the de Rham functor is essentially shown (compare the assertion (ii) of Theorem 5.4.1 written in \cite[p.~825]{KK1} with \cite[Proposition 3.3]{Mbid}). The essential surjectivity is not difficult to show by using the full faithfulness.
In \cite{Krh1}, \cite{Krh2}, a quasi-inverse is explicitly constructed.
\sk
(iv) In the algebraic case, the argument is more complicated, since we have the regularity {\it at infinity}, see \cite{Bo}, \cite{eq}, \cite{HT}, etc.
(This is essential to get a {\it canonical algebraic structure} on the vector bundle associated with a local system on a smooth complex variety by using the Deligne extension, see \cite{eq}.) The Riemann-Hilbert correspondence is used in an essential way in representation theory, see \cite{BB}, \cite{BK}. (This point does not seem to be sufficiently clarified in the last one.)
\msn
{\bf B.6.~$V$-filtration.} For a complex manifold $X$, set $Y:=X\times\C$ with $t$ the coordinate of $\C$.
We have the filtration $V$ on $\D_Y$ indexed by $\Z$ and such that
\msn
(i) $\,\,V^0\D_Y\subset\D_Y$ is the subring generated by $\Oc_Y$, $\dd_{y_i}$, and $t\dd_t$,
\msn
(ii) $\,\,V^j\D_Y=t^j\,V^0\D_Y$,\,\,\,$V^{-j}\D_Y=\msum_{0\les k\les j}\,\dd_t^k\,V^0\D_Y$\q($j\in\Z_{>0}$),
\msn
where the $y_i$ are local coordinates of $Y$, and $\dd_{y_i}:=\dd/\dd y_i$.
\sk
Let $M$ be a regular holonomic $\D_Y$-module. We say that $M$ is {\it quasi-unipotent} if so is $K:=\DR_Y(M)$, that is, if there is a stratification $\{S\}$ of $Y$ such that the restrictions of the cohomology sheaves $\Hc^iK$ to each stratum $S$ are $\C$-local systems having quasi-unipotent local monodromies around $\So\setminus S$.
\sk
For a quasi-unipotent regular holonomic left $\D_Y$-module $M$, there is a unique exhaustive filtration $V$ of Kashiwara \cite{Kvan} and Malgrange \cite{Ma3} indexed discretely by $\Q$ and satisfying the following three conditions:
\msn
(iii)\,\, $V^{\al}M$ ($\forall\,\al\in\Q$) are locally finitely generated $V^0\D_Y$-submodules,
\msn
(iv)\,\, $t\,V^{\al}M\subset V^{\al+1}M$ (with equality if $\al>0$),\,\,\,$\dd_t\,V^{\al}M\subset V^{\al-1}M\,\,$ for any $\al\in\Q$,
\msn
(v)\,\,\, $\dd_tt-\al$ is locally nilpotent on $\Gr_V^{\al}M\,\,$ $(\forall\,\al\in\Q)$.
\ms
Here we say that $V$ is {\it indexed discretely by} $\Q$ if there is a positive integer $m$ satisfying
$$V^{\al}M=V^{j/m}M\q\h{if}\q(j-1)/m<\al\les j/m\q\h{with}\q j\in\Z.
\leqno({\rm B}.6.1)$$
The existence of $V$ follows from that of $b$-functions in \cite{Khol2} where the holonomicity is actually sufficient (see \cite[Proposition 1.9]{rat}).
\sk
From condition (v) we can easily deduce the isomorphisms
$$t:\Gr_V^{\al}M\simto\Gr_V^{\al+1}M,\q\dd_t:\Gr_V^{\al+1}M\simto\Gr_V^{\al}M\q\q(\forall\,\al\ne 1),
\leqno({\rm B}.6.2)$$
It is also easy to show the following:
$$V^{\al}M=0\,\,\,(\forall\,\al>0)\q\h{if}\q{\rm supp}\,M\subset X\times\{0\}\subset Y,
\leqno({\rm B}.6.3)$$
$$t:V^{\al}M\simto V^{\al+1}M\q\q(\forall\,\al>0),
\leqno({\rm B}.6.4)$$
$$\h{$M\mapsto V^{\al}M$ (or $\Gr_V^{\al}M$) are exact functors ($\forall\,\al\in\Q$),}
\leqno({\rm B}.6.5)$$
see \cite[Lemma 3.1.3, Lemma 3.1.4, Corollary 3.1.5]{mhp}, where right $\D$-modules are used so that the action of $t\dd_t$ there corresponds to that of $-\dd_tt$ in this paper by (B.2.4), and an increasing filtration $V_{\ssb}$ is used there so that $V_{\al}=V^{-\al}$.
\sk
Set $Z:=g^{-1}(0)\subset X$. Let $M'_Z$ be the largest holonomic $\D_X$-submodule of $M$ supported in $Z$, and similarly for $M''_Z$ with submodule replaced with quotient module. Then we have the following canonical isomorphisms of $\D_X$-modules (see \cite[Proposition 3.1.8]{mhp}):
$$\aligned M'_Z&={\rm Ker}\bl(t:\Gr_V^0\Mt\to\Gr_V^1\Mt\br),\\ M''_Z&={\rm Coker}\bl(\dd_t:\Gr_V^1\Mt\to\Gr_V^0\Mt\br).\endaligned
\leqno({\rm B}.6.6)$$
Set $K:=\DR_Y(M)\in\HS(Y,\C)$. In the notation of (A.2), there are canonical isomorphisms
$$\aligned\DR_X(\Gr_V^{\al}M)&\simto\mh\psi_{t,{\mathbf e}(-\al)}K\q(\al\in(0,1]),\\ \DR_X(\Gr_V^0M)&\simto\mh\varphi_{t,1}K,\endaligned
\leqno({\rm B}.6.7)$$
such that $\exp(-2\pi i(\dd_tt-\al))$ on the left-hand side corresponds to the monodromy $T$ on the right-hand side, where ${\mathbf e}(-\al):=\exp(-2\pi i\al)$. In particular, $-(\dd_tt-\al)$ corresponds to $N:=(2\pi i)^{-1}\log T_u$ with $T=T_sT_u$ the Jordan decomposition. Moreover the morphisms
$$-\Gr_V\dd_t:\Gr_V^1M\to\Gr_V^0M,\q\Gr_Vt:\Gr_V^0M\to\Gr_V^1M
\leqno({\rm B}.6.8)$$
respectively correspond to the morphisms can and Var in (A.3.9) with $K'=K$ and $f=t$. (The sign before $\Gr_V\dd_t$ in (B.6.8) comes from the transformation between left and right $\D$-modules as in (B.2.4).)
\sk
In the case $M=(i_g)_*^{\D}\Oc_X$ with $i_g:X\into X\times\C$ the graph embedding by a holomorphic function $g$ on $X$, the proof of (B.6.7) is given in \cite{Ma3}, and this can be extended to the general regular holonomic case (see also the proof of \cite[Proposition 3.4.12]{mhp}).
There are canonical morphisms inducing the isomorphisms of (B.6.7) by using logarithmic functions, and these canonical morphisms are quite important for the proof of the stability theorem of Hodge modules by direct images.
\bs\bs
\vbox{\centerline{\bf Appendix C. Compatible filtrations}
\bsn
In this appendix we review some basics of compatible filtrations.}
\msn
{\bf C.1.~Compatible subobjects.} Let $\A$ be an abelian category, and $n\in\Z_{>0}$. We say that subobjects $B_i$ ($i\in[1,n]$) of $A\in\A$ are {\it compatible subobjects} if there is a {\it short exact $n$-ple complex} $K$ in $\A$ such that
\skn
(i) $\,\,K^p=0$ if $|p_i|>1$ for some $i\in[1,n]$, where $p=(p_1,\dots,p_n)\in\Z^n$.
\skn
(ii) $\,\,K^{p-\bo_i}\to K^p\to K^{p+\bo_i}$ is exact for any $p\in\Z^n$, $i\in[1,n]$.
\skn
(iii) $\,\,K^0=A$, $K^{-\bo_i}=B_i$ for any $i\in[1,n]$.
\msn
Here $\bo_i=((\bo_i)_1,\dots,(\bo_i)_n)$ with $(\bo_i)_j=\delta_{i,j}$
Note that conditions (i) and (ii) respectively correspond to ``short" and ``exact". In the case $n=2$, $K$ is the diagram of the {\it nine lemma}.
\msn
{\bf C.2.~Compatible filtrations.} We say that $n$ filtrations $F_{(i)}$ ($i\in[1,n]$) of $A\in\A$ form {\it compatible filtrations} if
$$F_{(1)}^{\nu_1}A,\,\dots\,,\,F_{(n)}^{\nu_n}A$$
are compatible subobjects of $A$ for any $\nu=(\nu_1,\dots,\nu_n)\in\Z^n$.
\sk
If $n=2$, then any $2$ filtrations $F_{(1)}$, $F_{(2)}$ are compatible filtrations. However, this does not necessarily hold for $n>2$.
\sk
We can show that if the $F_{(i)}$ ($i\in[1,n]$) form compatible filtrations of $A$, then their restrictions to $F^{\nu}A:=F_{(1)}^{\nu_1}\cdots F_{(n)}^{\nu_n}A$ also form compatible filtrations (see the proof of \cite[Corollary 1.2.13]{mhp}). Using the short exact complex $K$ with
$$K^0=F^{\nu}A,\q K^{-\bo_i}=F^{\nu+\bo_i}A\,\,\,\,(i\in[1,n]),$$
we can show that
$$\h{$\Gr_{F_{(1)}}^{\nu_1}\cdots\,\Gr_{F_{(n)}}^{\nu_n}A\,$ does not depend on the order of $\{1,\dots,n\}$.}
\leqno({\rm C}.2.1)$$
In fact, $\Gr_{F_{(i)}}^{\nu_i}$ corresponds to restricting $K$ to the subcomplex defined by $p_i=1$, and
$$\Gr_{F_{(1)}}^{\nu_1}\cdots\,\Gr_{F_{(n)}}^{\nu_n}A=K^{1,\dots,1}.
\leqno({\rm C}.2.2)$$
Note that (C.2.1) is not completely trivial even in the case $n=2$, where Zassenhaus lemma is usually used, and we can replace it by the diagram of the nine lemma as is explained above.
\msn
{\bf C.3.~Strict complexes.} Let $A^{\ssb}$ be a complex in $\A$ with $n$ filtrations $F_{(i)}$ ($i\in[1,n]$). We say that $\bl(A^{\ssb},F_{(i)}\,(i\in[1,n])\br)$ is {\it strict} if for any $j\in\Z$, $\nu=(\nu_1,\dots,\nu_n)\in\Z^n$, there is a short exact $n$-ple complex $K$ as in (C.1.1) such that
$$H^j\bl(\mcap_{i\in I}\,F_{(i)}^{\nu_i}A^{\ssb}\br)=K^{-\bo_I}\q\bl(\forall\,I\subset\{1,\dots,n\}\br),
\leqno({\rm C}.3.1)$$
where $\bo_I=((\bo_i)_I,\dots,(\bo_I)_n)$ with $(\bo_I)_j=1$ if $j\in I$, and $0$ otherwise (and $I$ can be empty in (C.3.1)).
\sk
We can show (see \cite[Corollary 1.2.13]{mhp}) that if $\bl(A^{\ssb},F_{(i)}\,(i\in[1,n])\br)$ is strict, then
$$\h{$H^j,\,\Gr_{F_{(1)}}^{\nu_1},\,\dots\,,\,\Gr_{F_{(n)}}^{\nu_n}$ on $\A^{\ssb}$ commute with each other $(\forall\,j\in\Z,\,\nu\in\Z^n)$.}
\leqno({\rm C}.3.2)$$
\vskip-7mm
$$\h{The induced filtrations $F_{(i)}$ ($i\in[1,n]$) on $H^jA^{\ssb}$ form compatible filtrations.}
\leqno({\rm C}.3.3)$$
\sk
Note that $\bl(A^{\ssb},F_{(i)}\,(i\in[1,2])\br)$ is strict if $(A^{\ssb},F_{(2)})$, $(\Gr_{F_{(2)}}^pA^{\ssb},F_{(1)})$ ($\forall\,p\in\Z)$ are strict, and $F_{(2)}^pA^{\ssb}=0$ for $p\gg 0$, where we assume that the filtered inductive limit in $\A$ is an exact functor, see \cite[Theorem~1.2.9]{mhp}.

\end{document}